\magnification=\magstep1


\font\typa=cmcsc10 scaled 833        
\font\typb=cmr10                     
\font\typc=cmbx10 scaled \magstep 1  
 2  
 3  

\font\eightrm=cmr8
\long\def\fussnote#1#2{{\baselineskip=9pt
     \setbox\strutbox=\hbox{\vrule height 7pt depth 2pt width 0pt}%
     \everypar{\hangindent=\parindent}
     \eightrm
     \footnote{#1}{#2}\everypar{}}}

\def\uw{{<{\kern-0.8em})}\ }
\def\ow{{<{\kern-0.5em})}\ }

\font\refp=cmr10
\font\refq=cmsl10
\def\ref#1 #2 #3 #4
{\par
\noindent
\refp#1\bu\
\refp#2:
\refq#3.
\refp#4.
\smallskip}

\def\qed{\vbox{\hrule
  \hbox{\vrule\hbox to 5pt{\vbox to
5pt{\vfil}\hfil}\vrule}\hrule}}

\def\cC{{\setbox0=\hbox{\rm C}\raise0.055\ht0\hbox to0pt
{\hskip0.3\wd0\vrule height0.85\ht0 width0.35pt\hss}\copy0}}

\def\cQ{{\setbox0=\hbox{\rm Q}\raise0.048\ht0\hbox to0pt
{\hskip0.29\wd0\vrule height0.88\ht0 width0.35pt\hss}\copy0}}

\def\mcc{{\setbox0=\hbox{$\scriptstyle{\rm C}$}\raise0.056\ht0\hbox to0pt
{\hskip0.28\wd0\vrule height0.76\ht0 width0.24pt\hss}\copy0}}

\def\mcq{{\setbox0=\hbox{$\scriptstyle{\rm Q}$}\raise0.051\ht0\hbox to0pt
{\hskip0.26\wd0\vrule height0.79\ht0 width0.24pt\hss}\copy0}}

\def\stroke#1{{\vrule height6.4pt width0.33pt} \kern-0.7pt{\rm #1}}
\def\Stroke#1{{\vrule height4.4pt width0.2pt} \kern-0.5pt{\scriptstyle{\rm #1}}}
\def\Strokeseven#1{{\vrule height4.4pt width0.2pt} \kern-0.5pt{{\rm #1}}}

\def\cR{\stroke{R}}

\def\cN{{\vrule height6.4pt width0.3pt} \kern-0.7pt{\rm N}}

\def\mcr{\Stroke{R}}

\def\mcn{\Stroke{N}}

\def\mch{{\vrule height4.4pt width0.19pt} \kern-0.4pt{\rm H}}

\def\cZ{{\raise6.5pt\hbox{{\vrule width5.6pt height0.35pt}}}
{\kern-6.7pt {\angle{\kern-5pt \angle}}}}

\def\mcZ{ \kern1.5pt{\hbox{\vrule width4.3pt
height0.22pt}{\raise1.4pt\hbox{${\kern-4.7pt
{\scriptscriptstyle{/}}{\kern-2.3pt{\scriptscriptstyle{/}}}}$}}
{\raise4.9pt\hbox{\kern-5.1pt{\vrule
width4.3pt height0.22pt}}}}\kern1pt}

\def\cU{{\setbox0=\hbox{\rm U}\raise0.055\ht0\hbox to0pt
{\hskip0.35\wd0\vrule height0.87\ht0
width0.35pt\hss}\copy0}}

\def\mcu{{\setbox0=\hbox{$\scriptstyle{\rm
U}$}\raise0.056\ht0\hbox to0pt {\hskip0.35\wd0\vrule
height0.89\ht0 width0.18pt\hss}\copy0}}

\def\C{\mathchoice{\cC}{\cC}{\mcc}{\mcc}}

\def\H{\mathchoice{\cH}{\cH}{\mch}{\mch}}

\let\N=\Nten

\let\R=\Rten

\def\Z{\mathchoice{\cZ}{\cZ}{\mcZ}{\mcZ}}

\hsize=159truemm
\vsize=245truemm

\def\3{\ss}
\def\\{\"}

\def\bb#1{{\raise1.5pt\hbox{$^=$}}\kern-6.5pt #1}

\def\h#1{{\hat{#1}}}

\def\rh#1{{\raise3.5pt\hbox{$\rightharpoonup$}}\kern-7.5pt #1}
\def\t#1{\tilde #1}


\def\ker{\mathop{\rm ker\,}\nolimits}

\def\log{\mathop{\rm log\,}\nolimits}

\def\rank{\mathop{\rm rank\,}\nolimits}

\font\dtsch=eufm10
\def\g{{\hbox{{\dtsch g}}}}
\def\J{{\hbox{{\dtsch J}}}}
\def\H{{\hbox{{\dtsch H}}}}
\def\V{{\hbox{{\dtsch V}}}}
\def\T{{\hbox{{\dtsch T}}}}
\def\t{{\hbox{{\dtsch t}}}}
\def\h{{\hbox{{\dtsch h}}}}

\def\( {\ )\kern-1,5pt(\kern-7pt-}
\def\doppelq{/ \hskip -2pt /}

\font\typc=cmbx10 scaled\magstep 1
\font\typb=cmr10 scaled\magstep 1
\font\typa=cmr10 scaled\magstep 3
\font\meinfoot=cmsl8
\font\meinft=cmr8

\vsize=8.3 true in

\font\dtsch=eufm10
\def\g{{\hbox{{\dtsch g}}}}
\def\J{{\hbox{{\dtsch J}}}}
\def\H{{\hbox{{\dtsch H}}}}
\def\V{{\hbox{{\dtsch V}}}}
\def\T{{\hbox{{\dtsch T}}}}
\def\t{{\hbox{{\dtsch t}}}}
\def\h{{\hbox{{\dtsch h}}}}

\def\( {\ )\kern-1,5pt(\kern-7pt-}
\def\doppelq{/ \hskip -2pt /}

\font\typc=cmbx10 scaled\magstep 1
\font\typb=cmr10 scaled\magstep 1
\font\typa=cmr10 scaled\magstep 3
\font\meinfoot=cmsl8
\font\meinft=cmr8

\vsize=8.3 true in

{\centerline {\typa Isoparametric submanifolds and}}

\bigskip{\centerline{\typa a Chevalley-type restriction theorem}}

\vskip 1,5 cm
{\centerline{\typb Ernst Heintze, Xiaobo Liu, Carlos
Olmos\footnote{$^1$}{\meinft Research supported by Universidad Nacional de C\'ordoba,
partially supported by CONICOR, Secyt-UNC and CIEM}
}
{\footnote{}
{\meinfoot 1991 Mathematics Subject Classification. 53B25, 53C30, 53C40\hfil
\break
\meinft Key words and phrases. Isoparametric submanifolds, Chevalley restriction
theorem, polar actions, Riemannian submersions, representations of compact Lie
groups, submanifolds of Hilbert space\hfil}
}

\vskip 1 cm

\par{\narrower\noindent
{\meinft
ABSTRACT. We define and study isoparametric submanifolds of general ambient
spaces and of arbitrary codimension. In particular we study their behaviour
with respect to Riemannian submersions and their lift into a Hilbert space.
These results are used to prove a Chevalley type restriction theorem which
relates by restriction eigenfunctions of the Laplacian on a compact Riemannian
manifold which contains an isoparametric submanifold with flat sections to
eigenfunctions of the Laplacian of a section. A simple example of such an
isoparametric foliation is given by the conjugacy classes of a compact Lie
group and in that case the restriction theorem is a (well known) fundamental result in
representation theory. As an application of the restriction theorem we show
that isoparametric submanifolds with flat sections in compact symmetric spaces
are level sets of eigenfunctions of the Laplacian and are hence related to
representation theory. In addition we also get the following results.
Isoparametric submanifolds in Hilbert space have globally flat normal bundle,
and a general result about Riemannian submersions which says that focal
distances do not change if a submanifold of the base is lifted to the total
space.
}\par}

\vskip 2 cm
\noindent
{\typc 1. Introduction}

\bigskip
A {\sl hypersurface} $M$ of a Riemannian manifold is called {\sl isoparametric}, if
its parallel manifolds have constant mean curvature or equivalently if $M$ is
given (locally) as the regular level set of a function $f$ whose two {\sl
parameters} $\Vert grad f\Vert$ and  $\Delta f$ are functions of $f$. The
parallel manifolds of $M$ are also isoparametric and one gets (locally) two
orthogonal foliations: The isoparametric foliation which consists of $M$ and
its parallel manifolds and the $1$-dimensional foliation by the geodesics
perpendicular to $M$. Typical examples are the principal orbits of
cohomogeneity$-1$ actions. In spaces of constant curvature the definition
implies that even the principal curvatures of $M$ are
constant and not only their sum (E. Cartan). But of course this is no longer
true in more general ambient manifolds.

The aim of this paper is twofold. Our first goal is to propose a definition for
isoparametric submanifolds of arbitrary codimension which extends the above
definition for hypersurfaces in a natural way. We call a submanifold $M$ to be {\sl isoparametric},
if $\nu M$ is flat, the (close
by) parallel manifolds have constant mean curvature in radial directions and
if $M$ admits sections, that is, if for each $p\in M$ there exists a totally
geodesic submanifold which is tangent to $\nu_pM$.
Again one gets locally two transversal foliations, namely the isoparametric
foliation by $M$ and its parallel manifolds and orthogonally the totally
geodesic foliation by the sections.
We study some properties of
these submanifolds, in particular their behaviour with respect to Riemannian
submersions. We also show that our definition coincides with that of Terng
([T1], [T2]), if the ambient manifold is a euclidean space, possibly an
infinite dimensional Hilbert space. Furthermore we prove that a submanifold of
a compact symmetric space is isoparametric with flat sections
if and only if it is equifocal in the sense of Terng and
Thorbergsson [TT].

Simple examples of isoparametric submanifolds are the (non-totally
geodesic)\break
leaves of a warped product. More interesting, the principal orbits of a {\sl
polar} action are isoparametric, where we call an isometric action to be polar, if
the distribution of the normal spaces to the principal orbits is integrable.
The integral manifolds of this distribution are necessarily totally geodesic and yield the
sections. A typical example of a polar action is a compact Lie group $G$ with
biinvariant metric which acts on itself by conjugation. The conjugacy classes
of regular elements are the isoparametric submanifolds and the maximal tori
are the sections in this case.

The linearized version of this example, namely the adjoint action of $G$ on
its Lie algebra (which is also polar), brings us to our second topic, the
{\sl Chevalley-type restriction theorem} for isoparametric submanifolds. If $\g$
denotes the Lie algebra of $G\ ,\ \t\subset\g$ a maximal abelian subalgebra
and $W$ the Weyl group, then one version of the classical Chevalley
restriction theorem says (cf. [Hel2], Ch II, \S 5), that restricting polynomials from $\g$ to $\t$
defines an isomorphism from the space of $AdG$-invariant polynomials on $\g$
onto the space
of $W$-invariant polynomials on $\t$ or
equivalently from the space of polynomials on $\g$ which are constant on the
leaves of the isoparametric foliation onto the space of polynomials on $\t$
which are constant on the intersections of the leaves with $\t$. Now let $X$
be a compact symmetric or more generally normal homogeneous space and
$M\subset X$ a compact isoparametric submanifold with flat sections whose
parallel manifolds decompose $X$ (which is always the case if $X$ is simply
connected). If $\Sigma$ is a section, then we have the following
Chevalley-type restriction theorem (Theorem 7.1).

\bigskip\noindent
{\bf Theorem} {\sl Restricting functions from $X$ to $\Sigma$ defines an
isomorphism from $C^\Delta(X)^M$ onto $C^\Delta(\Sigma)^M$.}

\bigskip
Here $C^\Delta(X)^M$ denotes the space of finite sums of eigenfunctions of the
Laplacian of $X$ which are constant on $M$ and its parallel manifolds.
Similarly $C^\Delta(\Sigma)^M$ denotes the space of finite sums of
eigenfunctions of the Laplacian on $\Sigma$ which are constant on the
intersections of $\Sigma$ with the parallel manifolds of $M$. The finite sums
of eigenfunctions replace the polynomials in the classical restriction
theorem. In fact, if one replaces $\g$ and $\t$ in this theorem by the
corresponding unit spheres $S(\g)$ and $S(\t)$ - which are invariant under $G$
and $W$, respectively - then one  gets the equivalent theorem that
restriction defines an isomorphism from $C^\Delta(S(\g))^G$ onto
$C^\Delta(S(\t))^W$. Recall that the finite sums of eigenfunctions of the
Laplacian on the unit sphere are precisely the restrictions of the polynomials
of the ambient euclidean space.

\bigskip
Even in the example above - the action of $G$ on itself by conjugation - the
theorem is interesting, yet well known. In fact, it gives a weak form of
the fundamental result of representation theory which says that restriction
defines an isomorphism from $R(G)$ onto $R(T)^W$ where $T$ is a maximal torus
and $R(G)$ and $R(T)$ denote the representation rings of $G$ and $T$,
respectively. More precisely, it gives an isomorphism between the corresponding
vector spaces $R(G)\otimes\C$ and $(R(T)\otimes\C)^W$ (recall that the
representation rings are $\Z$-modules).

\bigskip
If the ambient space
$X$ in the theorem is in addition simply connected, then it follows moreover
that $C^\Delta(X)^M$ is a polynomial algebra in $k=codim M$ generators and that
these generators can be chosen as eigenfunctions of the Laplacian of $X$
(Theorem 7.6). This in turn implies that isoparametric submanifolds with flat
sections in compact, simply connected, normal homogeneous
 spaces are the simultanous level sets of eigenfunctions of
the Laplacian (Corollary 7.7) and are hence related to representations of the
isometry group of $X$. Recall, that the pull back of an eigenfunction of
$X=G/H$ to $G$ is the coefficient of some representation of $G$.
This result might be even interesting for isoparametric hypersurfaces of compact
symmetric spaces and explains to some extent why all the known examples of isoparametric
hypersurfaces in spheres are connected with representations, namely
Clifford representations and  isotropy representations of symmetric spaces.

\bigskip
One of our principal tools for the proof of the above results is the lifting
of an isoparametric submanifold with flat sections in certain spaces to an
isoparametric submanifold of a Hilbert space. This idea goes back to Terng and
Thorbergsson [TT]. For the price of infinite dimensions, the ambient manifold
is so to speak linearized and the theory of isoparametric submanifolds of
Hilbert spaces (as developed by Terng in [T2]) becomes applicable. In particular
one gets a formula for the mean curvature of the original isoparametric
submanifold in this way ([KT]), which does not seem to be obtainable so far in
a direct manner. This formula
plays a key role for our results, as it relates the
Laplacian of a section $\Sigma$ to the Laplacian of the ambient manifold $X$.
We develop these ideas of Terng - Thorbergsson and King - Terng carefully in
Chapter 5 and 6 and put them into a broader context. We also
generalize a result of Terng - Thorbergsson ([TT], Lemma 5.12) and O'Neill
([O'N2], Theorem 4) to the effect that for any Riemannian submersion $\pi:\hat X\to X$ focal
distances and multiplicities of a submanifold  $M$ of $X$
coincide with that of its lift $\hat M=\pi^{-1}(M)\subset\hat X$.
 This gives for example a simple geometric proof for the fact that fibres
of a Riemannian submersion  $\pi:V\to X$ are minimal
(in a certain regularized sense), if $V$ is a Hilbert and $X$ a symmetric
space.

To prove the equivalence of our notion of isoparametric submanifold with that
of Terng for submanifolds of Hilbert spaces we essentially have to extend
Cartan's theorem  that isoparametric submanifolds have
constant principal curvatures to infinite dimensions (Theorem 4.2). But there is also a more
technical point, as Terng requires isoparametric submanifolds of Hilbert
spaces to have {\sl globally}
flat normal bundle, that is with trivial, not only finite holonomy, and our
definition is purely local. However, we show in Appendix B that global
flatness is a consequence of her other conditions, like in finite dimensions.

\bigskip
Finally, we prove in Appendix A that our definition of polar actions (the
distribution of normal spaces to the principal orbits is integrable) implies the
existence of totally geodesic complete immersed sections through each point which
meet all orbits and always perpendicularly. But in general these sections will
neither be closed nor embedded (as required in the definition of [PT1]).
Counterexamples are already given by certain cohomogeneity$-1$ actions on
compact, simply connected symmetric spaces.

\vskip 1 cm

\noindent
{\typc 2. Isoparametric submanifolds}

\bigskip
Let $M$ denote an immersed submanifold of a Riemannian
manifold $N$.
$M$ is said to have globally flat normal bundle, if the holonomy of the
normal bundle is trivial, i.e. if any normal vector can be extended to a
globally defined parallel normal field. It then follows in particular that $\nu M$ is
flat. On the other hand, if $\nu M$ is flat, then there exists for each
$p\in M$ an open neighborhood $U$ of $p$ in $M$, over which the normal bundle is
globally flat. By eventually restricting $U$ further one can find an $r>0$
such that in addition  the exponential map
 is a diffeomorphism on $\nu^r M_{\vert U}$, where
$\nu^r M=\{\xi\in \nu_ M\mid \Vert\xi\Vert<r\}$. The "local tube"
$V:=\exp(\nu^rM_{\vert U})$ around $M$ is then foliated by the parallel
submanifolds $U_\xi=\{\exp\xi(q)\vert q\in U\},$ where $\xi$ is a parallel
normal field along $U$ of length less than $r$.
 On $V\setminus U$ there exists the radial vector field
${\partial\over\partial r}=\hbox{grad} r$, where $r$ is the distance to $M$.
It follows from the first variation formula that the restriction of
$\partial\over\partial r$ to any $U_\xi$ is a normal vector field. We say
that {\sl locally the parallel submanifolds of $M$ have constant mean
curvature in radial directions} if for each $p\in M$ we can find $U$ and $r$
as above, such that the mean curvature in the direction of $\partial\over
\partial r$, i.e. $tr\ A_{\partial\over\partial r}$, is constant along
$U_\xi$ for all parallel normal fields along $U$ with $\Vert\xi\Vert\in(0,r)$,
where $A$ denotes the shape operator of $U_\xi$.

\bigskip\noindent
{\bf Definition} $M$ is called {\sl almost isoparametric}, if $\nu M$ is
flat and if locally the parallel submanifolds of $M$ have constant mean
curvature in radial directions.

\bigskip\noindent
{\bf Remark:} It follows by continuity that an almost isoparametric
submanifold has constant mean curvature in the direction of any parallel
normal field, and thus has parallel mean curvature.

\bigskip
The condition on the mean curvature of parallel submanifolds can also be
described in the following way.

\bigskip\noindent
{\bf Proposition 2.1} {\sl Let $M$ be an immersed  submanifold with flat
normal bundle. Then $M$ is almost isoparametric if and only if for each $p\in
M$ there exists a local tube $V=\exp(\nu^rM\vert_U)$ around $M$ such that
for all parallel $\xi$ in $\nu^r M_{\vert U}$ the mapping\quad $U\to U_\xi\
,\ x\mapsto\exp\xi(x),$ is volume preserving up to a constant factor.}

\bigskip\noindent
{\bf Proof}: Let $\xi_0:=\xi(x_0)$ for some $x_0\in U$ and
$Y_1,\dots,Y_m\ (m=\dim M)$ a basis of $M$-Jacobi fields along $c_{\xi_0}$,
i.e. of Jacobi fields associated to variations of $c_{\xi_0}$ through
geodesics starting orthogonally from $M$, where $c_v$ is the geodesic in the
direction of $v$. Then
$$
(\log\Vert Y_1\wedge\dots\wedge Y_m\Vert)'(t)=-tr\ A_{{\partial\over\partial
r}(c_{\xi_0}(t))}\ ,\leqno(*)
$$
as one can choose the $Y_i$ in such a way that they are orthonormal at
$t$ which implies
$$\eqalign{
&(\log\Vert Y_1\wedge\dots\wedge Y_m\Vert)'(t)={1\over 2}(\log\langle
Y_1\wedge\dots\wedge Y_m,Y_1\wedge\dots\wedge Y_m\rangle)'(t)=\cr
&=\langle(Y_1\wedge\dots\wedge Y_m)'\ ,\ Y_1\wedge\dots\wedge Y_m\rangle
(t)={\sum\limits_i}\langle Y_i',Y_i\rangle(t)=\cr
&=-{\sum\limits_i}\langle A_{\partial\over\partial r}Y_i,Y_i\rangle\ .\cr}
$$
The last equality comes from the fact, that geodesics starting orthogonally
from $U$ hit $U_{t\xi}$ orthogonally as well, so that the $Y_i$ are also
$U_{t\xi}$-Jacobi fields.
Since the ratio of volume elements of $U_{t\xi_0}$
at $c_{\xi_0}(t)$ and $U$ at $x_0$ is given by the ratio of $\Vert
Y_1\wedge\dots\wedge Y_m\Vert$ at $t$ and $0$, respectively, the proposition
follows from $(*)$. \hfill\qed

\bigskip
In general the parallel submanifolds of an almost isoparametric submanifold
need not be almost isoparametric themselfes. To remedy this situation we
will require in addition that $M$ has totally geodesic sections.

\bigskip\noindent
{\bf Definition} An immersed submanifold $M$ is said to have {\sl
(totally geodesic) sections} if for each $p \in M$ there exists a totally
geodesic submanifold $\Sigma$ (also called a section) which meets $M$ at $p$
orthogonally and whose dimension is equal to the codimension of $M$.

\bigskip
Obviously such a section has to coincide locally with $\exp(\nu_pM).$
Therefore $M$ has totally geodesic sections if and only if  for each $p\in
M\ ,\ \exp(\nu_pM)$ is totally geodesic in a neighborhood of $p$. This
condition is automatically satisfied if $M$ is a hypersurface or if $N$ has
constant sectional curvature.

The following special case will be of interest later.

\bigskip\noindent
{\bf Definition} $M$ is said to have {\sl flat sections} if
$M$ has totally geodesic sections
which are flat with respect to the induced metric.

\bigskip
Besides hypersurfaces and submanifolds of euclidean spaces, interesting
examples of submanifolds with flat sections are principal orbits of
hyperpolar actions (cf. [HPTT] and the remark below).

\bigskip\noindent
{\bf  Definition} An immersed submanifold is called {\sl isoparametric}
if it is almost isoparametric and if it has totally geodesic sections.

\bigskip
The distinction between isoparametric and almost isoparametric is
redundant if $M$ is a hypersurface or if the ambient manifold has constant
curvature. Important examples of isoparametric submanifolds arise as the
principal orbits of {\sl polar} actions. We call (following a suggestion of W. Ziller)
an isometric action of a Lie
group on a Riemannian manifold to be polar if the distribution of the normal
spaces to the principal orbits is integrable. The integral manifolds are then
totally geodesic and yield the sections. If these are flat with respect to the
induced metric, the action is called hyperpolar ([HPTT]). Recently A. Kollross
[Ko] has classified hyperpolar actions on irreducible compact symmetric spaces
and Podest\`a and Thorbergsson [PTh] have classified polar actions on compact
rank $1$ symmetric spaces. Our definition of polar is somewhat weaker than
that of Palais and Terng [PT1] in that we do not require the existence of a
closed embedded section. However, we prove in Appendix A that also with our
definition sections (in complete manifolds) can be always extended to complete
totally geodesic immersions which meet all orbits and always perpendicularly.
But in general they are neither closed nor embedded. This shows for example
the diagonal action of $SO(3)$ on $S_{r_1}^2\times S^2_{r_2}$ which is of
cohomogeneity one and thus polar. But the geodesics normal to the principal
orbits are not closed if the ratio $r_{1/\displaystyle{r_2}}$ of the radii of the spheres
is not rational.

\bigskip
The following two results show the strong influence of the existence of
totally geodesic sections for submanifolds with flat normal bundle.

\bigskip\noindent
{\bf Proposition 2.2} {\sl Let $M$ be a submanifold with globally flat
normal bundle and assume that there exists $r>0$ such that
$\exp_{\vert\nu^rM}$
is a diffeomorphism and $\exp(\nu^r_pM)$ is totally geodesic for all $p\in
M$. Then the parallel manifolds $M_\xi,\Vert\xi\Vert<r$, meet the sections
$\exp(\nu^r_pM)$ orthogonally as well and have globally flat normal
bundle, too. In fact, a parallel normal field along $M$ transported to
$M_\xi$ by parallel translation along the geodesics $\exp t\xi$ is a
parallel normal field along $M_\xi$.}

\bigskip\noindent
{\bf Remark}: It follows in particular that  the
radial vector field  ${\partial\over\partial r}$ is parallel in the normal
bundle of any $M_\xi,0<\Vert\xi\Vert<r$.

\bigskip\noindent
{\bf Proof}:

(i) Let $\xi$ be a parallel normal field along $M$ with
$\Vert\xi\Vert<r$ and fix some $p\in M$. Let $c(t):=\exp t\xi(p)$ and
$q:=c(1)\in M_\xi$. If $v\in T_pM$ and $\gamma$ is a curve in $M$ with
$\dot\gamma(0)=v$, then $\alpha(s,t):=\exp t\xi(\gamma(s))$ is a
variation
of
$c$ through geodesics with $\alpha(s,0)=\gamma(s)$ and $\alpha(s,1)\in
M_\xi$. The variation vector field $Y_v(t):={\partial\alpha\over\partial
s}(0,t)$ is the Jacobi field along $c$ with $Y_v(0)=v$ and $Y'_v(0)={D\over
\partial t}{\partial\alpha\over\partial s}(0,0)={D\over\partial
s}{\partial\alpha\over\partial t}(0,0)=\xi'(0)=-A_\xi v\in T_pM$. Since
$T_qM_\xi=\{Y_v(1)\vert v\in T_pM\},$ we have to show that these Jacobi
fields
$Y_v$ are perpendicular to the section $\Sigma:=\exp \nu_pM$. But this
follows from the fact that $Y_v$ is a solution of $Y''+R_{\dot c}\ Y=0$ with
initial conditions in $(T_{\dot c(0)}\Sigma)^\bot\ ,$ where $R_{\dot
c}=R(.,\dot c)\dot c\ ,$ and that $T_{\dot c(0)}\Sigma$ as well as $(T_{\dot
c(t)}\Sigma)^\bot$ are families of parallel subspaces along $c(t)$ which are
invariant under $R_{\dot c}$. Therefore $M_\xi$ meets $\Sigma$
perpendicularly.

(ii) Let $\eta$ be another parallel normal field along $M$ and
$\eta^*$
be the vector field along $M_\xi$ which is obtained by parallel translating
$\eta(p)$ along $\exp t\xi(p)$ to $\exp\xi(p)$. Then $\eta^*$ is a normal
field along $M$ by the above considerations. Let $p\in M$ and $\gamma(s)$ a
curve in $M$ starting at $p$ and
$\alpha(s,t)=\exp t\xi(\gamma(s))$. Let
$\Sigma=\exp_p\nu_pM$ as above and $\eta(s,t)$ be the parallel translation
of $\eta(\gamma(s))$ along $\tau\to\exp\tau\xi(\gamma(s))$ up to $\exp t
\xi(\gamma(s))$, so that $\eta(s,0)=\eta(\gamma(s))$ and
$\eta(s,1)=\eta^*(\exp\xi{\gamma(s))}$. If $v(t)$ is parallel along $c(t)$
and tangent to $\Sigma$, then ${\partial\over\partial
t}\langle{D\eta\over\partial s},v\rangle_{/_{s=0}}=\langle{D\over\partial
t}{D\eta\over\partial s},v\rangle_{/_{s=0}}=\langle R(\dot
c,{\partial\alpha\over\partial s})\eta,v\rangle_{/_{s=0}}=0$, since
${\partial\alpha\over\partial s}_{/_{s=0}}\ \bot\ \Sigma$ and $\dot
c,\eta,v$ are tangent to $\Sigma$. Since $\langle{D\eta\over\partial
s},v\rangle_{/_{s=0}}=0$ at $t=0\ ,\ {D\eta\over\partial s}_{/_{s=0}}$ is
perpendicular to $\Sigma$ for all $t\in[0,1]$, in particular for $t=1,$
which shows that $\eta^*$ parallel in $\nu M_\xi$. \hfill\qed

\bigskip\noindent
{\bf Proposition 2.3} {\sl The following conditions on an immersed
submanifold $M$ of $N$ are equi\-valent.

\item{(i)} $M$ has flat normal bundle and totally geodesic sections.

\item{(ii)} Locally, (i.e. for each $p\in M$ there exists a neighborhood in
$N$, such that in this neighborhood) $M$ is the fibre of a Riemannian
submersion $f:N\to B$ with integrable horizontal spaces.

\item{(iii)} Locally, $M$ is the leaf of a foliation which admits an
orthogonal, transversal and totally geodesic foliation.

\item{(iv)} Locally, $N$ splits as $N=N_1\times M$ with metric $g_1\oplus
g_2(x)$, where $g_1$ is a fixed metric
\item{} on $N_1$ and $g_2(x)$ is a metric on $M$ which depends on $x\in N_1$.}

\bigskip\noindent
{\bf Proof}: "(i) $\Rightarrow$ (iii)" follows from Proposition 2.2.,
"(iv) $\Rightarrow$ (ii)" is obvious, as the projection onto $N_1$ is a
Riemannian submersion. "(ii) $\Rightarrow$ (i)" follows from the remark that
the O'Neill tensor of $f:N\to B$ vanishes, and Lemma 3.1 below. Thus we are
left to show "(iii) $\Rightarrow$ (iv)". Since the foliations are
transversal, $N$ is locally diffeomorphic to $N_1\times M$ and the $N_1$
factors are by assumption totally geodesic and perpendicular to the
$M$-factors. If $c_1$ and $c_2$ are differentiable curves in $N_1$ and $M$,
respectively and $\alpha(s,t):=(c_1(s),c_2(t))\ ,$ we therefore have
${1\over
2}{\partial\over\partial t}\Vert{\partial\alpha\over\partial
s}\Vert^2=\langle{D\over\partial t}{\partial\alpha\over\partial s},{\partial
\alpha\over \partial s}\rangle=\langle{D\over\partial
s}{\partial\alpha\over\partial t},{\partial\alpha\over\partial
s}\rangle={\partial\over\partial s}\langle{\partial\alpha\over\partial
t},{\partial\alpha\over\partial
s}\rangle-\langle{\partial\alpha\over\partial t},{D\over\partial
s}{\partial\alpha\over\partial s}\rangle=0$. Hence the metric on the
$N_1$-factors is constant. \hfill\qed

\bigskip
The following theorem shows that for an isoparametric submanifold the mean
curvature of the parallel submanifolds is not only constant in radial
directions but also in the direction of any parallel normal field. More
precisely we have:

\bigskip\noindent
{\bf Theorem 2.4} {\sl Let $M$  be an immersed submanifold  of $N$ with
flat normal bundle and totally geodesic sections. Then $M$ is isoparametric
if and only if locally the parallel submanifolds have parallel mean
curvature.}

\bigskip\noindent
{\bf Proof}: We may assume that $M$ has globally flat normal bundle
and there exists $r>0$ such that $\exp$ is a diffeomorphism on $\nu^rM$ and
$\exp\nu_p^rM$ is totally geodesic for all $p\in M$. If the close by
parallel manifolds $M$ have parallel mean curvature, then they have
constant mean curvature in the radial direction $\partial\over\partial r$,
as $\partial\over\partial r$ is parallel in $\nu M_\xi$ by the remark
following Proposition 2.2. Hence $M$ is isoparametric.

On the other hand, let $M$ be isoparametric and let $\eta$ be the vector
field on $\exp\nu^r M$ whose restriction to any parallel manifold $M_\xi$
is the mean curvature vector field of $M_\xi$. By Proposition 2.4 we may
assume that $\exp\nu^rM=N_1\times M$ with metric $g_1\oplus g_2(x),x\in
N_1$, and that $N_1=B_r(0)\subset \nu_pM$. We then have to show that
$\eta=\eta(x,y)$, where $x\in N_1$ and $y\in M$, is constant in $y$. By
eventually restricting $M$ further we may assume that there exist on $M$
globally defined vector fields $X_1,\dots,X_m$ which form at every point a
basis of the tangent space. We view these also as vector fields on
$N_1\times M$ which are tangent to $M$ and constant in the $N_1$-direction.
For any $\xi\in T_{(x,y)}N_1\times M$ which is normal to $M$ we then have
as in the proof of Proposition 2.1
$$
tr\ A_\xi=-\xi\log\Vert X_1\wedge\dots\wedge X_m\Vert\ ,
$$
where $A_\xi$ is the shape operator of $M$. Hence $\eta=\hbox{grad}_{N_1}
f,$
where $f=-\log\Vert X_1\wedge\dots\wedge X_m\Vert$ and where the subscript
$N_1$ denotes the orthogonal projection onto $TN_1$. If $X$ and $Y$ are
vector fields on $N_1\times M$ with $X$ tangent to $M$ and constant in the
$N_1$ direction and $Y$ tangent to $N_1$ and constant in the $M$-direction,
then $[X,Y]=0$ and $X\langle\eta,Y\rangle=X(Yf)=Y(Xf)$. If $Y$ is the
radial vector field, then $\langle\eta,Y\rangle=tr\ A_y$ is constant in the
direction of $X$ by assumption. Hence $Xf$ is constant in radial directions
and thus constant on any section $N_1\times\{q\},q\in M$. This implies
$Y(Xf)=0$ for any $Y$ as above and hence $X\langle\eta,Y\rangle=0$. Thus
$\eta=\eta(x,y)$ is constant in $y$. \hfill\qed

\bigskip
From Proposition 2.3 and the theorem we get:

\bigskip\noindent
{\bf Corollary 2.5} {\sl The locally defined parallel manifolds of an
isoparametric submanifold are isoparametric as well, and thus define locally
an isoparametric foliation.}

\bigskip
Together with Proposition 2.1 we get:

\bigskip\noindent
{\bf Corollary 2.6} {\sl An immersed submanifold $M$ is isoparametric if
and only if $\nu M$ is flat, $M$ has totally geodesic sections and if the
projection  from  $M$ to any (sufficiently close) parallel
submanifold along the sections is volume preserving up to a constant
factor.}

\bigskip
In particular the leaves of a warped product $N\times_f M\ ,\ f:N\to\R$,
more precisely the factors $\{n\}\times M$ are isoparametric, as the projection from
one such leaf to another is even an isometry up to a constant factor
and the factors $N\times\{ m\}$ are totally geodesic.

\bigskip
We finally discuss the relation of our definition to some of the previous
ones. The classical definition used by Segre, Levi-Civita and E. Cartan (cf.
[Ca]) calls a {\sl hypersurface} $M$ of a {\sl space form} $N$ to
be isoparametric if it is given as a regular level set of a function
$f:N\to\R$ whose two
"parameters" $\Vert\hbox{grad}f\Vert$ and $\Delta f$ are functions of $f$.
Geometrically this means precisely that parallel submanifolds have constant
mean curvature (and that the normal bundle is globally flat). Thus the
definition is essentially
equivalent to ours. Cartan had shown that a
hypersurface of a space form is isoparametric if and only if it has constant
principal curvatures. Terng [T1] used this result to define a submanifold of
arbitrary codimension of a {\sl space form} to be isoparametric if it has
flat
normal bundle and if the principal curvatures in the direction of any
(locally defined) parallel normal field are constant. It follows from her
extension of Cartan's result above, that for submanifolds of {\sl space
forms}
this definition is equivalent to ours. She also extended  her definition
to submanifolds in Hilbert space in [T2] and we will prove (Corollary 4.3)
that also in this case both definitions coincide, by extending Cartan's
result further. On the other hand, Wang had already shown in [W] that an
isoparametric hypersurface in ${\C}P^n$ does not need to have constant principal
curvatures. He constructs an example by pushing down a certain
non-homogeneous isoparametric hypersurface in $S^{2n+1}$ via the Hopf
fibration into $\C P^n$. This
shows also that the property of having constant principal curvatures for a
hypersurface is not preserved under Riemannian submersions.
Another definition of isoparametric
submanifolds, but with the
restriction $\dim M\ge{1\over 2}\dim N$ and $\nu M$ globally flat, had been
given by Carter and West [CW] using differential forms and the Hodge
star-operator. Their Theorems
(1.1) to (2.5) show that it is equivalent to ours under the restrictions made.
Finally we  mention that already Harle had given a definition of
isoparametric submanifold in [Ha] which is very similar to ours, but had
applied it only to submanifolds of space forms. We would like to thank G.
Thorbergsson for pointing out this reference to us.

\vskip 1 cm

\noindent
{\typc 3. Isoparametric submanifolds and Riemannian submersions}

\bigskip
In later sections  we will make strong use of the fact, that isoparametric
submanifolds behave nicely with respect to  Riemannian submersions. Let
$\pi:X\to B$ be a Riemannian submersion.
The tangent bundle of $X$ splits
orthogonally into the vertical bundle $\cal V$ (which is tangent to the
fibres)
and the horizontal bundle ${\cal H}={\cal V}^\bot$. Let $\nabla$ and
$\bar\nabla$ be the covariant derivatives of $X$ and $B$, respectively. Then
we have for all vector fields $V,W$ on $B$
$$
\nabla_{\hat V}\hat W=\bar\nabla_VW+O'N(\hat V,\hat W)\ ,
$$
where $\hat{}$ denotes the horizontal lift and
$O'N(\hat V,\hat W):=(\nabla_{\hat V}\hat W)^{\cal V}$ is one of the tensors
of O'Neill ([O'N]), usually denoted by A. Since also $O'N(\hat V,\hat W)={1\over 2}
[\hat V,\hat W]^{\cal V}$, this tensor vanishes precisely, if the horizontal
distribution is
integrabel and hence totally geodesic. The curvature $\bar
K(\sigma)$ of a plane $\sigma$ tangent to $B$ is related to the curvature
$K(\hat\sigma)$ of a horizontal lift by
$$
\bar K(\sigma)=K(\hat\sigma)+3\Vert O'N(\hat v,\hat w)\Vert^2\ ,
$$
where $v,w$ is an orthonormal base of $\sigma$. Thus these two curvatures
coincide if and only if $O'N$ vanishes on $\hat \sigma$.

The vanishing of $O'N$ can also be described as follows.

\bigskip\noindent
{\bf Lemma 3.1} {\sl Let $b\in B$. Then $O'N$ vanishes along the fibre
$F_b:=\pi^{-1}(b)$ if and only if for each $v\in T_bB$ the horizontal lift
$\hat v$ of $v$ is parallel in $\nu F_b$.}

\bigskip\noindent
{\bf Proof}: We extend $v$ locally to a vector field which we also denote
by $v$. Then $\hat v$ is parallel in $\nu F_b$ if and only if
$\langle\nabla_Y\hat v,\hat w\rangle=0$ for all vertical vector fields $Y$
and all horizontal lifts $\hat w$ of vectors $w\in T_bB$. Since $Y$ is
$\pi$-related to zero and $\hat v$ to $v\ ,\ [Y,\hat v]$ is $\pi$-related to
zero, i.e. vertical. Thus $\langle\nabla_Y\hat v,\hat
w\rangle=\langle\nabla_{\hat v}Y,\hat w\rangle=-\langle Y,\nabla_{\hat
v}\hat w\rangle=-\langle Y,O'N(\hat v,\hat w)\rangle$, from which the lemma
follows. \hfill\qed

The following Proposition generalizes a result of Wang ([W], Corollary 2) from
hypersurfaces to submanifolds of arbitrary codimension.

\bigskip\noindent
{\bf Proposition 3.2} {\sl Let $\pi:X\to B$ be a Riemannian submersion with
minimal fibres and $M\subset B$ an embedded submanifold. Let
$M^*:=\pi^{-1}(M)$ and assume that $O'N=0$ on $\nu M^*$. Then $M$ is almost
isoparametric if and only if $M^*$ is almost isoparametric.}

\bigskip\noindent
{\bf Remarks}:

\noindent(i) Since $M^*$ contains with any point the entire fibre through
the point, its normal spaces are horizontal and are mapped isometrically by
$\pi_*$ onto the normal spaces of $M$. In particular $M^*$ has the same
codimension as $M$.

\noindent
(ii) If $M$ is a hypersurface, so is $M^*$ and the condition
on $O'N$ is automatically satisfied as O'N $(X,X)=0$ for all $X$.

\bigskip\noindent
{\bf Proof}: It follows from $O'N=0$ on $\nu M^*$ and the proof of Lemma
3.1 that the horizontal lift $\hat \xi$ of a normal vector field $\xi$ is
parallel in $\nu M^*$ in vertical directions. Since $(\nabla_{\hat
v}\hat \xi)^\bot=\widehat{(\bar\nabla_v\xi)^\bot}$ for any $v\in TB\ ,\
\hat\xi$ is also parallel in $\nu M^*$ in horizontal directions if and only
if $\xi$ is parallel in $\nu M$.
Hence $\nu M^*$ is flat if and only if $\nu M$ is flat and in that case the
horizontal lifts of parallel normal fields are precisely the parallel
normally fields of $M^*$ locally.

Let $\nu M$ be flat and $V=\exp(\nu^rM_{\vert U})$ a local tube around $M$.
Let
$\xi\in \nu^r M_{\vert U}$ be parallel, $U^*=\pi^{-1}(U)$ and
$U^*_\xi:=\pi^{-1}(U_\xi)$. Since $\pi(\exp t\hat\xi)=\exp t\xi$, the
parallel manifold $(U^*)_{\hat\xi}$ of $U^*$ is contained in $U^*_\xi$. But
actually equality $(U^*)_{\hat\xi}=U^*_\xi$ holds, as one can reverse
the roles of $U$ and $U_\xi$ for this argument: $U^*$ is obtained from
$(U^*)_{\hat\xi}$ by lifting the geodesics $\exp t\xi$ from $U$ to $U_\xi$
with the reversed direction horizontally. The proposition follows therefore
from the following lemma. \hfill\qed

\bigskip\noindent
{\bf Lemma 3.3} {\sl Let $\pi:X\to B$ be a Riemannian submersion with minimal
fibres, $M\subset B$ be an embedded submanifold and $M^*:=\pi^{-1}(M)$. Then
the mean curvature vector field of $M^*$ is the horizontal lift of that of
$M$. In particular $tr\ A^*_{\hat \xi}=tr\ A_\xi$ for any normal vector
$\xi$ of $M$, where $\hat\xi$ is a horizontal lift and $A,A^*$ denote the
shape operators of $M$ and $M^*$, respectively.}

\bigskip\noindent
{\bf Proof}: Let $\hat b\in X$ and $b:=\pi(\hat b)$. We choose an
orthonormal basis $e_1,\dots,e_m\in T_bM$ and extend their horizontal lifts
$\hat e_1,\dots,\hat e_m$ by an orthonormal basis $f_1,\dots,f_k$ of
$T_{\hat b}(\pi^{-1}(b))$ to
one of $T_{\hat b}M^*$. Then $\eta^*={\sum\limits_i}\alpha^*(\hat
e_i,\hat
e_i)+{\sum\limits_j}\alpha^*(f_j,f_j)={\sum\limits_i}\alpha^*(e_i,e_i)$,
where $\eta^*$ denotes the mean curvature vector field and $\alpha^*$ the
second fundamental form of $M^*$. Since $\langle\alpha^*(\hat e_i,\hat
e_i),\hat \xi\rangle=\langle\nabla_{\hat e_i}\hat e_i,\hat
\xi\rangle=\langle\bar\nabla_{e_i}e_i,\xi\rangle=\langle\alpha(e_i,e_i),\xi\rangle$
for any $\xi\in\nu_bM$ the lemma
\break follows. \hfill\qed

\bigskip\noindent
{\bf Remark}: If $\nu M$ is globally flat in Proposition 3.2, so is $\nu
M^*$, since the horizontal lift of a parallel normal field is a parallel
normal field of $M^*$. The converse is not always true, but is true, if the
fibres are connected. In that case the restriction of a parallel normal
field $\xi^*$ in $\nu M^*$ to a fibre is not only locally a horizontal lift
of a  normal vector of $M$, but also globally, and $\xi^*$ can be pushed
down to a globally defined parallel normal field of $M$.

\bigskip\noindent
{\bf Theorem 3.4} {\sl Let $\pi:X\to B$ be a Riemannian submersion with
minimal fibres, $M\subset B$ be a submanifold and $M^*:=\pi^{-1}(M)$. Then
the following conditions are equivalent.}

(i)$M^*$ is isoparametric with horizontal sections, i.e. sections
which are perpendicular to the fibres.

(ii) $M$ is isoparametric and $O'N=0$ on all horizontal lifts of
tangent vectors to sections of $M$.

\bigskip\noindent
{\bf Proof}: A horizontal section $\Sigma^*$ of $M^*$ is mapped under $\pi$
isometrically onto a section of $M$, as geodesics in $\Sigma^*$ are
horizontal and are thus mapped to geodesics in $M$. From O'Neill's
curvature formula it thus follows that $O'N=0$ on tangent vectors to
$\Sigma^*$. Together with Proposition 3.2 this proves "(i) $\Rightarrow$
(ii)".

If (ii) holds then $M^*$ is almost isoparametric and its parallel manifolds
are horizontal lifts of parallel manifolds of $M$ (locally) by Proposition
3.2 and its proof. The normal distribution $\cal D$ to the parallel
manifolds of $M^*$ is therefore horizontal. It is also integrabel and
totally geodesic,
since we have for any $\xi,\eta\in{\cal D}_p\quad(\nabla_\xi\eta)^{\cal
V}=0$ by the assumption on $O'N$ and $\langle\nabla_\xi\eta,\hat
v\rangle=\langle\bar\nabla_{\pi_*\xi}\pi_*\eta,v\rangle=0$ for any $v$ which
is tangent to the parallel submanifold of $M$ through $\pi(p)$. Thus (i)
follows. \hfill\qed

\bigskip\noindent
{\bf Corollary 3.5} {\sl Let $\pi:X\to B$ be a Riemannian submersion with
minimal fibres and let $X$ have non negative curvature. If $M\subset
B$ is an isoparametric submanifold with flat sections, so is
$M^*=\pi^{-1}(M)$.}

\bigskip\noindent
{\bf Proof}: By O'Neill's curvature formula, condition (ii) of the theorem
is satisfied. \hfill\qed

\bigskip\noindent
{\bf Remark}: The Corollary as well as Proposition 3.2 and Theorem 3.4 are
also true if $X$ is a separable Hilbert space. In fact, we will show in
section 5 that the mean curvature of $M^*$ can be defined also in this
situation and in such a way that Lemma 3.3 is valid. But then the same
arguments apply.

\bigskip
In a special situation  also the converse to Corollary 3.5 holds.

\bigskip\noindent
{\bf Corollary 3.6} {\sl Let $G$ be a Lie group with biinvariant metric,
$H\subset G$ a closed subgroup and $G/H$ the normal homogeneous space
endowed
with the quotient metric, so that $\pi:G\to G/H$ is a Riemannian submersion.
Then a submanifold $M$ of $G/H$ is isoparametric with flat sections if and
only if $M^*=\pi^{-1}(M)$ is isoparametric with flat sections.}

\bigskip\noindent
{\bf Proof}: Let $M^*$ be isoparametric with flat sections and assume for
convenience that $e\in M^*$. Then a section at $e$ is an abelian Lie group
$A$ which is perpendicular to the fibre $H$. Since the metric is invariant
under left translations, $A$ is also perpendicular to the fibres $a\cdot H\
,\ a\in A$. Hence the first condition of the theorem is satisfied and $M$ is
isoparametric with flat sections.

The other direction follows from Corollary 3.5. \hfill\qed

\vskip 1 cm
\goodbreak

\noindent
{\typc 4. Isoparametric submanifolds of Hilbert space}

\bigskip
In this section we extend our definition of isoparametric submanifolds
to the situation where the ambient
manifold is a  Hilbert space and show that it coincides with
that of Terng in [T2]. This may be viewed as a generalization of theorems of
Cartan and Terng to the infinite dimensional case.

We will always (in the whole paper) denote by $V$ a real separable Hilbert
space. We begin by defining the regularized trace of certain compact operators,
which will later be the shape operators of submanifolds of $V$. Our
regularization is somewhat different from
the $\zeta$-regularization used by King and Terng [KT] and
has the advantage of being easier to handle. In all
relevant cases, however, both regularizations will coincide (cf. the remark below).

\bigskip\noindent
{\bf Definition} Let $A:V\to V$ be a compact selfadjoint operator with non-zero
eigenvalues
$\mu_1\le\mu_2\le\dots<0<\dots\le\lambda_2\le\lambda_1$, repeated with multiplicities.
Then $A$ is
called {\sl regularizable} if $tr\ A^2<\infty$ and
${\sum\limits^\infty_1}(\lambda_k+\mu_k)$ converges where $\lambda_k$ or
$\mu_k$ is understood to be zero, if there are less than $k$ positive or
negative eigenvalues, respectively. If $A$ is regularizable, then
$$
tr_rA:={\sum\limits^\infty_1}(\lambda_k+\mu_k)
$$
is called the {\sl regularized trace} of $A$.

\bigskip
In general $\Sigma\lambda_k$ and $\Sigma\mu_k$ will not converge. A typical
example is $\lambda_k=-\mu_k=1/k$ with $tr_rA=0$. The definition of the
regularized trace therefore depends strongly on the chosen order of the
eigenvalues.

\bigskip\noindent
{\bf Remark}: King and Terng call a compact, self-adjoint operator $A:V\to V$
with non-zero eigenvalues
$\mu_1\le\mu_2\le\dots<0<\dots\le\lambda_2\le\lambda_1\quad\zeta$-regularizable
if ${\sum\limits_k}\lambda_k^s-{\sum\limits_k}\vert\mu_k\vert^s$ converges
for all $s>1$ and $tr_\zeta A:={\lim\limits_{s\searrow
1}}({\sum\limits_k}\lambda_k^s-{\sum\limits_k}\vert\mu_k\vert^s)$ exists. By the
dominated convergence theorem, it follows easily that $tr_\zeta A=tr_rA$, if $A$
is $\zeta$-regularizable and ${\sum\limits_k}(\lambda_k+\mu_k)$ converges absolutely.
In fact, ${\lim\limits_{s\searrow
1}}{\sum\limits_k}(\lambda^s_k-\vert\mu_k\vert^s)={\sum\limits_k}{\lim\limits_{s\searrow
1}}(\lambda_k^s-\vert\mu_k\vert^s)={\sum\limits_k}(\lambda_k+\mu_k)$, as $\vert
x^s-y^s\vert\le 2\vert x-y\vert$ for all $s\in[1,2]$ and $x,y\in[0,1]$.

\bigskip
\noindent{\bf Definition} An immersed submanifold $M$ of $V$ is called {\sl
regularizable}, if for each $\xi\in\nu M$ the shape operator $A_\xi$ is
regularizable. In this case $tr_rA_\xi$ is also called the mean curvature of
$M$ in direction of $\xi$.

\bigskip
The following lemma implies that  parallel manifolds of $M$ are regularizable
in radial directions, if $M$ is regularizable. It also motivates the condition
$tr\ A^2<\infty$ in the definition of regularizability.

\bigskip
\noindent{\bf Lemma 4.1} {\sl Let $\{a_k\}_{k\in\N}$ be a sequence of complex
numbers such that $\Sigma a_k$ and $\Sigma\vert a_k\vert^2$ converge. Then
also $\Sigma{a_k\over 1-a_kz}$ converges for all $z\in\C$ for which $a_kz\ne
1$ for all $k\in\N$. The limit function $f$ is meromorphic in $\C$ and
$\{1/a_k\mid k\in\N\ ,\ a_k\ne 0\}$ is the set of its poles. The residue of
$f$ at a pole $1/a_k$ is equal to $-\sharp\{\ell\in\N\mid a_\ell=a_k\}$. In
particular, $f$ determines the non-zero $a_k's$ together with their
multiplicities.}

\bigskip\noindent
{\bf Proof}: Given $R>0$ there exists $k_0\in\N$ with $\vert
a_k\vert\le{1\over 2R}$ for all $k\ge k_0$ and hence $\vert
1-a_kz\vert\ge{1\over 2}$ in the closed ball of radius $R$ around zero.
Therefore ${\sum\limits^\infty_{k_0}}{a_k\over
1-a_kz}={\sum\limits^\infty_{k_0}}a_k+{\sum\limits^\infty_{k_0}}{a^2_kz\over
1-a_kz}$ converges in this ball uniformely to a holomorphic function. The
other statements are now obvious.\hfill\qed

\bigskip
If $M$ is a regularizable submanifold of a Hilbert space $V$ and $\xi$ is a
parallel normal field along $M$, then, as in finite dimensions,
$M_{t\xi}:=\{p+t\xi(p)\mid p\in M\}$ is, at least locally and for small $t$,
a submanifold of $V$ diffeomorphic to $M$. In fact, $p\mapsto p+t\xi(p)$ has
differential $id-tA_{\xi(p)}$ and is thus for any $p\in M$ injective in some
neighborhood of $p$, if $t$ is small. Moreover the tangent spaces of $M$ and
$M_{t\xi}$ coincide in corresponding points so that $\xi^*_t$, which is
defined to be the parallel translation of $\xi(p)$ to $p+t\xi(p)$, is a
parallel normal field along $M_{t\xi}$. Again as in finite dimensions the
eigenvalues of the shape operator of $M_{t\xi}$ in direction of $\xi^*_t$ are
$\lambda_k\over 1-t\lambda_k$ and $\mu_k\over 1-t\mu_k$, respectively, where
$\mu_1\le\mu_2\le\dots\le 0\le \lambda_2\le\lambda_1$ are the eigenvalues of
$A_\xi$ and $t$ is sufficiently small. But then ${\mu_1\over
1-t\mu_1}\le{\mu_2\over 1-t\mu_2}\le\dots\le 0\le\dots\le{\lambda_2\over
1-t\lambda_2}\le{\lambda_1\over 1-t\lambda_1}$ due to the monotonicity of the
function $x\mapsto{x\over 1-tx}$. Hence the regularized mean curvature of
$M_{t\xi}$ in direction of $\xi^*_t$ is defined by Lemma 4.1. Moreover, if
this is constant along $M_{t\xi}$ for small $t$, then Lemma 4.1 shows that the
non-zero eigenvalues of $A_\xi$ together with their multiplicities are
constant along $M$.

\bigskip\noindent
{\bf Definition} Let $M$ be an immersed submanifold of a separable Hilbert
space. $M$ is called {\sl isoparametric}, if $M$ is regularizable, $\nu M$ is
flat and if for each $p\in M$ there exists a neighborhood in which the
parallel manifolds $M_\xi$ have constant mean curvature in direction
$\xi^*$ (which is the parallel translation of $\xi$ to $M_\xi$ and thus a
constant multiple of the radial direction).

\bigskip
The following result generalizes the characterization of Cartan and Terng of
iso\-parametric submanifolds in euclidean space to infinite dimensions.

\bigskip\noindent
{\bf Theorem 4.2} {\sl Let $M$ be an immersed regularizable submanifold of a
separable Hilbert space with flat normal bundle. Then the following conditions
are equivalent:

(i) $M$ is isoparametric.

(ii) If $\xi$ is a parallel normal field in an open connected subset
$U$ of $M$ then the shape operators $A_{\xi(p)}$ and $A_{\xi(q)}$ are
orthogonally equivalent for all $p,q,\in U$.}

\bigskip\noindent
{\bf Proof}: Let (iii) be the statement: If $\xi$ is a parallel normal field
in an open connected subset $U$ of $M$ then the non-zero eigenvalues
of $A_\xi$ together with their multiplicities are constant on $U$.

Lemma 4.1 together with the discussion above proves (i) $\Rightarrow$ (iii).
Vice versa, (i) is a trivial consequence of (iii) since the mean curvature of
$M_\xi$ in the direction of $\xi^*$ is given by
${\sum\limits_k}({\lambda_k\over 1-\lambda_k}+{\mu_k\over 1-\mu_k})$ where
$\mu_1\le\mu_2<0<\dots\le\lambda_2\le\lambda_1$ are the non-zero
eigenvalues of $A_\xi$. It is also obvious that (ii) implies (iii).

Under the assumption of (iii) it follows from perturbation theory that the
eigen\-space to a non-zero eigenvalue varies continuously with $p\in M$. In
fact, the orthogonal projection onto such an eigenspace can be described as a
certain integral over the resolvent of $A_{\xi(p)}$ (cf. Kato [K], p. 178 -
181) and thus depends smoothly on $p$. Hence also the zero-eigenspace as the
orthogonal complement of the others varies continuously and is therefore
either infinite  dimensional or has finite constant dimensions. This proves
that (i) follows from (iii).\hfill\qed

\bigskip
We are now in a position to compare our definition of isoparametric submanifold
with that of Terng [T2]. Terng calls an immersed submanifold $M$ in a
separable Hilbert space $V$ {\sl proper Fredholm}, if the endpoint map
$Y:=\exp_{\nu M}:M\to V$ is Fredholm (i.e. has a Fredholm derivative) and the
restriction of $Y$ to each normal disk bundle of finite radius $r$ is proper.
While the first condition is a purely local one, and is in fact equivalent to
the compactness of the shape operators of $M$, the second condition is a global one.
It implies in particular that the immersion itself is proper and hence that
$M$ is complete as a metric space. Now, Terng defines an immersed, proper
Fredholm submanifold of a Hilbert space of finite codimension to be
isoparametric, if $\nu M$ is globally flat (that is has trivial holonomy, not only finite) and
if for any parallel normal field $\xi$ the shape operators $A_{\xi(p)}$ and
$A_{\xi(q)}$ are orthogonally equivalent for all $p,q\in M$.

\bigskip\noindent
{\bf Corollary 4.3} {\sl Let $M$ be an immersed proper Fredholm submanifold of
finite codimension of a separable Hilbert space. Then $M$ is isoparametric
(according to our definition) if and only if $M$ is isoparametric in the sense
of Terng.}

\bigskip\noindent
{\bf Proof}: In Appendix B we show that {\sl global}
flatness of the normal bundle in Terng's definition can be replaced by
flatness of $\nu M$.
Hence, if $M$ is isoparametric, then it is also isoparametric  in the
sense of Terng by Theorem 4.2. The other direction also follows immediately
from the Theorem except that the regularizability of $M$ is not clear a
priori.
 But this is a consequence of the explicit description of the eigenvalues
of $A_\xi$ given by Terng [T2], cf. also Lemma 7.5 and its proof.

\bigskip
For later applications we note the following two results about the regularized
trace of compact, self-adjoint operators.

\bigskip\noindent
{\bf Lemma 4.4} {\sl Let $V_i\ ,\ i=1,\dots,k$, be separable Hilbert spaces and
$A_i:V_i\to V_i$ compact, self adjoint, regularizable operators. If
$V:={\bigoplus\limits^k_1}V_i$ and $A:={\bigoplus\limits^k_1}A_i$, then $A:V\to V$ is
regularizable and $tr_rA={\sum\limits^k_1}tr_rA_i$.}

\bigskip\noindent
{\bf Proof}: By induction, it suffices to consider the case $k=2$. Let
$\mu_1\le\mu_2\le\dots<0<\dots\le\lambda_2\le\lambda_1$ and
$\mu_1'\le\mu_2'\le\dots<0<\dots\lambda_2'\le\lambda_1'$ be the non-zero
eigenvalues, repeated with multiplicities, of $A_1$ and $A_2$ respectively. For
any $N\in\N$ we choose $r,r',s$ and $s'$, such that $r+r'=s+s'=N$ and
$\lambda_1,\dots,\lambda_r,\lambda_1',\dots,\lambda_r'$ are the $N$ largest
and $\mu_1,\dots,\mu_s,\mu'_1,\dots,\mu_{s'}$ the $N$ smallest eigenvalues of
$A$. If none of the sequences
$\{\mu_i\},\{\mu'_i\},\{\lambda_i\},\{\lambda_i'\}$ is finite, then
necessarily $r,r',s$ and $s'$ tend to infinity together with $N$. From
$$
{\sum\limits^s_1}(\lambda_i+\mu_i)+{\sum\limits^{s'}_1}(\lambda'_i+\mu_i')
\le{\sum\limits^r_1}\lambda_i+{\sum\limits^{r'}_1}\lambda_i'+{\sum\limits^s_1}\mu_i+
{\sum\limits^{s'}_1}\mu_i'\le{\sum\limits^r_1}(\lambda_i+\mu_i)+{\sum\limits^{r'}_1}
(\lambda_i'+\mu'_i)
$$
the lemma follows, as the left and right hand side of the inequality both
tend to $tr_rA_1+tr_rA_2$. If e.g. $\{\mu_i\}$ is finite, say
$\mu_1\le\mu_2\le\dots\le\mu_{s_0}<0$, and the other sequences are not, then $r,r'$ and
$s'$ tend to infinity, $\Sigma\lambda_i$ converges and $r\cdot\lambda_r$ tends
to zero. From $s'=N-s_0\ge N-r=r'$ for large $N$ and
${\sum\limits^r_1}\lambda_i+{\sum\limits^{r'}_1}\lambda'_i+{\sum\limits^{s_0}_1}\mu_i+
{\sum\limits^{s'}_1}\mu_i'={\sum\limits^{s_0}_1}(\lambda_i+\mu_i)+{\sum\limits^r_{s_0+1}}\lambda_i+
{\sum\limits^{s'}_1}(\lambda_i'+\mu_i')-{\sum\limits^{s'}_{r'+1}}\lambda_i'$
the lemma follows also in this case, as
${\sum\limits^{s'}_{r'+1}}\lambda_i'\le(s'-r')\lambda'_{r'+1}=(r-s_0)\lambda'_{r'+1}\le(r-s_0)\lambda_r$
and hence ${\sum\limits^{s'}_{r'+1}}\lambda'_i$ tends to zero. If more than
one of the sequences is finite, the lemma is rather obvious.\hfill\qed

\bigskip
The final result is essentially due to King and Terng ([KT], Theorem 4.2) who
proved it for the $\zeta$-regularization. In fact, our proof is just a
simplified version of theirs, adapted to our regularization.

\bigskip\noindent
{\bf Lemma 4.5 (King and Terng)} {\sl Let $V$ be a separable Hilbert space and
$A,B:V\to V$ self-adjoint operators with $A$ compact, regularizable and $B$ of
finite rank. Then $A+B$ is regularizable and $tr_r(A+B)=tr_rA+tr\ B$.}

\bigskip\noindent
{\bf Proof}: By induction, we may assume that rank $B=1$ and hence
$B(x)=\epsilon<x,v>v$ for some $v\in V$ and $\epsilon\in\{\pm 1\}$. Let
$v={\sum\limits_i}v_i$ be a decomposition of $v$ into eigenvectors of $A$ to
pairwise different eigenvalues. If this sum is finite, the lemma follows
immediately, as $V$ then splits orthogonally into $A-$ and $B-$ invariant
subspaces $V_0\oplus V_1$ with $\dim V_0<\infty$ and $B(V_1)=0$. If the sum
${\sum\limits^\infty_{i=1}}v_i$ is infinite, we let
$w_N:={\sum\limits^N_{i=1}}v_i$ and $B_N(x)=\epsilon{<x,w_N>}w_N$. Then
$tr_r(A+B_N)=tr_rA+tr\ B_N$ and $tr\ B_N=\Vert w_N\Vert^2$ converges to $tr\
B=\Vert v\Vert^2$ (Note that Lemma 4.4 (b) of [KT], which is used there to
prove ${\lim\limits_{N\to\infty}}tr\ B_N=tr\ B$, is not correct). Hence we are
left to show that $A+B$ is regularizable and
${\lim\limits_{N\to\infty}}(tr_rA+B_N)=tr_r(A+B)$.

Let $\mu_1\le\mu_2\le\dots<0<\dots\le\lambda_2\le\lambda_1$ and
$\mu_1(N)\le\mu_2(N)\le\dots<0<\dots\le\lambda_2(N)\le\lambda_1(N)$ be the
non-zero eigenvalues of $A+B$ and $A+B_N$, respectively. Then
$\mu_{i-1}\le\mu_i(A+B_N)\le\mu_{i+1}$ and
$\lambda_{i+1}\le\lambda_i(A+B_N)\le\lambda_{i-1}$ for all $N$ and all $i\ge
2$. In fact, this follows from
$\lambda_{n+m+1}(T_1+T_2)\le\lambda_{n+1}(T_1)+\lambda_{m+1}(T_2)$ for any
compact, self-adjoint operators $T_1,T_2$ ([DS]) and the fact, that $A+B_N$
differs from $A+B$ only by an operator of rank $1$ (cf. Lemma 4.5 of [KT]).
Hence
$$
\lambda_l-\lambda_{k+1}+\mu_k-\mu_{l+1}\le{\sum\limits^l_{i=k}}(\lambda_i+
\mu_i)-{\sum\limits^l_{i=k}}(\lambda_i(N)+\mu_i(N))\le\lambda_k-\lambda_{l+1}+\mu_l-\mu_{k-1}
$$
for all $N$ and all $l\ge k\ge 2$. This proves that $A+B$ is regularizable and
that for all $\epsilon>0$ there exists $k_0$ such that
$\vert{\sum\limits^\infty_{k_0}}(\lambda_i+\mu_i)-{\sum\limits^\infty_{k_0}}(\lambda_i(N)+\mu_i(N))\vert
<\epsilon/2$ for all $N$. Since $\lambda_i(N)\to\lambda_i$ and
$\mu_i(N)\to\mu_i$ due to the convergence of $A+B_N$ to $A+B$ in norm (cf.
[DS]), we finally get $\vert
tr_r(A+B)-tr_r(A+B_N)\vert=\vert{\sum\limits^\infty_1}(\lambda_i+\mu_i)-
{\sum\limits^\infty_1}(\lambda_i(N)-\mu_i(N))\vert$ for all large
$N$.\hfill\qed

\vskip 1 cm

\noindent
{\typc 5. Lifting isoparametric submanifolds to a Hilbert space}

\bigskip
Isoparametric submanifolds with flat sections of finite dimensional manifolds
$X$ can be  studied for certain, but important $X$, by
lifting them to a Hilbert space and exploiting the theory of isoparametric
submanifolds there. This idea, which will be of decisive importance for our
restriction theorem, goes back to Terng and Thorbergsson [TT] in the
case of equifocal submanifolds. The goal of this chapter is  to describe and
simplify their results and to put them into a broader context.

The essential assumption on the finite dimensional ambient manifold $X$
 will be the following: $X$ is the base of a Riemannian
submersion $\pi:V\to X$ with minimal fibres, where $V$ denotes (as always) a real
separable Hilbert space.
The most important example of such a manifold $X$ is a Lie group $G$ with
biinvariant metric. In fact, let $V_G$ be the set of $H^1$-curves\quad $g:[0,1]\to G$ with
$g(0)=e$. Then $V_G$ may be identified with the Hilbert space $L^2([0,1],\g)$
of $L^2$-curves in the Lie algebra $\g$ of $G$ by mapping $g$ onto
$u(t):=g'(t)g^{-1}(t)$, i.e. onto the curve of tangent vectors pulled back by
right translations to the identity. If one considers $V_G$ in this way as a
Hilbert space (as we will do), then the endpoint map $\pi_G:V_G\to
G,g\mapsto g(1)$, becomes a Riemannian submersion ([TT]). Moreover, $\pi_G$ has minimal
fibres according to King and Terng ([KT]). In the next chapter we will give a simple
geometric explanation for this fact.

This example can be extended as follows.
If $H\subset G$ is a closed subgroup then the canonical projection $G\to G/H$
becomes in a natural way a Riemannian submersion with minimal (actually
totally geodesic) fibres, too. Therefore also the composition $\pi:V_G\to G\to G/H$
has this property (the minimality of the fibres follows from
Lemma 5.2 below). In particular any compact symmetric
space and more generally any normal homogeneous space $X$ is the base of a
Riemannian submersion $\pi:V\to X$ with minimal fibres, where $V$ is a Hilbert
space.

If one considers instead of $G/H$ more generally double quotients $G\doppelq H$ where
$H\subset G\times G$ operates freely on $G$ by $(h_1,h_2).
g=(h_1gh_2^{-1})$ then again $G\to G\doppelq H$ and hence also the composition
$V_G\to G\doppelq H$ is a Riemannian submersion with respect to the natural metric on $G\doppelq H$.
But its fibres are in general not minimal.

All of the examples above of Riemannian submersions $\pi:V\to X$  are of the form
$V\to V/P$ where $P$ is a group of isometries acting freely on
$V$ with orbits of finite codimension. The double quotients $G\doppelq H$ arise for example by taking
$P:=P(G,H)$ to be the group of $H^1$-curves \quad $h:[0,1]\to G$ with $(h(0),h(1))\in
H$, which acts isometrically on $V_G$ by $h.g(t)=h(t)g(t)h^{-1}(0)$. This group
has been introduced and studied by Terng in [T3]. We do not know whether there
are further examples. In this respect it might be interesting to note that
Gromoll and Walchap have recently classified  Riemannian submersions $\pi:V\to
X$ with $\dim V<\infty$, that is with $V=\R^n$ ([GW1],[GW2]). They show that
these are precisely of the above form $V\to V/P$, where $P$ is a group of
glide transformations. Note also that the base $X$ of a Riemannian submersion
$\pi:V\to X$, where $V$ is a Hilbert space,
 has non-negative curvature by O'Neill's formula and that most
known examples of manifolds with non-negative curvature arise in this way (at
least up to diffeomorphism).

Let $\pi:V\to X$ be a Riemannian submersion onto a finite dimensional
Riemannian manifold $X$. If $f:M\to X$ is an isometric immersion we
denote by $\hat f:\hat M\to V$ its lift to $V$ where $\hat M=\{(p,v)\in
M\times V\mid f(p)=\pi(v)\}$ and $\hat f(p,v)=v$. If $M$ is embedded then
$\hat M$ can be of course identified with $\pi^{-1}(M)$. Note that $\hat f$ is
an immersion with the same codimension as $f$ and that the normal spaces
of $\hat f$ are horizontal. We endow $\hat M$ with the metric induced by $\hat
f$ and observe that the natural projection $\hat M\to M$ then becomes a
Riemannian submersion, too. We recall that a differentiable mapping between
Hilbert manifolds is called Fredholm if its differentials are Fredholm and
that an immersed submanifold $N$ of $V$ is called proper Fredholm if $\exp_{\nu N}:\nu N\to V$
is Fredholm and the restriction of
$\exp_{\nu N}$ to normal vectors of bounded length is proper. It has been shown
by Terng in [T2] that $\exp_{\nu N}$ is Fredholm iff the shape operators of $N$
are compact.

\bigskip\noindent{\bf Lemma 5.1} {\sl Let $\pi:V\to X,f:M\to X$ and $\hat f:\hat M\to
V$ be as above. Then:

(i) The normal exponential map $\exp_{\nu\hat M}:\nu\hat M\to V$ is Fredholm.

(ii) $\hat M$ is proper Fredholm iff $f:M\to X$ is proper.}

\bigskip\noindent{\bf Proof}:

(i) If $c:[0,1]\to X$ is a geodesic segment then the horizontal lifts
of $c$ connect the fibres $\pi^{-1}(c(0))$ and $\pi^{-1}(c(1))$ with each
other, and therefore induce a diffeomorphism between these fibres. It follows
that for any $\hat\xi_0\in\nu\hat M\ ,\ \exp_{\nu\hat M}:\nu\hat M\to V$ maps the
submanifold $\{\hat \xi\in\nu\hat M\mid\pi_*\hat\xi=\pi_*\hat\xi_0\}$
diffeomorphically onto the fibre through $\exp\hat\xi_0$. Since this
submanifold and the fibre both have finite codimension in $\nu \hat M$ and
$\hat M$, respectively, the differential of
$\exp_{\nu\hat M}$ at $\hat\xi_0$ and hence $\exp_{\nu\hat M}$ itself are
Fredholm.

(ii) We only prove that "$f$ proper" implies "$\hat M$ proper Fredholm",
the other direction being almost trivial. Because of (i) it thus suffices to
show that any sequence $\hat\xi_n\in\nu\hat M$ with $\Vert\hat \xi_n\Vert\le
r$ and $\exp\hat \xi_n\to v$ has a convergent subsequence. Recall that
$\hat\xi_n$ is a pair $(\hat p_n,\bar v_n)$ where $\hat p_n=(p_n,u_n)\in\hat M$
and $\bar v_n\in T_{\hat f(\hat p_n)}V$, and thus $\bar v_n=(u_n,w_n)$. Let
$\Phi_t:TV\to TV$ be the geodesic flow of $V$, i.e. $\Phi_t(u,w)=(u+tw,w)$. Then
the vectors $\Phi_1(\bar v_n)$ are horizontal and of length $\le r$. By
assumption their foot points $\exp\bar v_n=\exp\hat\xi_n$ converge to $v$.
Hence a subsequence $\Phi_1(\bar v_{n_k})$ converges, say to $(v,w)$. Therefore
$\bar v_{n_k}$ converges to $\Phi_{-1}(v,w)$. Since $f(p_{n_k})=\pi(u_{n_k})$
converge and $f$ is proper we also may assume that the $p_{n_k}$ converge
which finally proves the convergence of the $\hat\xi_{n_k}$. \hfill\qed

\bigskip\noindent{\bf Remark}: A different  proof of (i) will follow from Corollary
6.2, which relates the eigenvalues of the shape operators of $\hat M$ to the focal
points of $M$.

\bigskip
The following result is essentially due to King and Terng who proved it in
case $X=G$ a compact Lie group and $V=V_G$ ([KT], Theorems 4.12 and 4.17).

\bigskip\noindent{\bf Lemma 5.2} {\sl Let $\pi:V\to X,f:M\to X$ and $\hat f:\hat M\to
V$ be as above, and assume in addition that $\pi$ has minimal fibres. Then
$\hat M$ is regularizable. Moreover the regularized mean curvature vector field of $\hat
M$ is defined and coincides with the horizontal lift of that of $M$.}

\bigskip\noindent{\bf Proof}: Since this result is local we may assume that $M$
is embedded. Let $\hat p\in\hat M$ and $\hat\xi\in \nu_{\hat p}\hat M$. Then
$T_{\hat p}\hat M$ may be identified with $T_{\hat p}F\oplus T_pM$ where $F$
is the fibre through $\hat p$ and $p=\pi\hat p$. The shape operator $A^{\hat
M}_{\hat\xi}$ of $\hat M$ therefore decomposes into $A+B+C$ where $A$ and $B$
preserve the decomposition of $T_{\hat p}\hat M$ with $A(T_pM)=B(T_{\hat
p}\hat F)=0$, while $C$ interchanges the summands. We may identify $A$ with
the shape operator $A^F_{\hat\xi}$ of the fibre and $B$ with the shape
operator $A^M_\xi$ of $M$ in direction $\xi:=\pi_*\hat\xi$. Since the rank of
$B+C$ is at most $2\cdot\dim T_pM$ and hence finite, Lemma 4.5  implies that
$A^{\hat M}_{\hat\xi}$ is regularizable and
$$
tr_r A^{\hat M}_{\hat \xi}=tr_r A+tr(B+C)=tr_r A^F_{\hat \xi}+tr A^M_\xi+tr
C=tr A^M_\xi\ .
$$
Therefore $\hat\xi\to tr_rA^M_{\hat\xi}$ is linear and we may define the mean
curvature vector field $\hat \eta$ of $\hat M$ by $<\hat\eta,\hat\xi>:= tr_r
A^{\hat M}_{\hat\xi}$. Since this latter expression coincides with $tr
A^M_\xi=<\eta,\xi>$ where $\eta$ denotes the mean curvature vector field of
$M\ ,\ \hat\eta$ is the horizontal lift of $\eta$.\hfill\qed

\bigskip
As always, we assume that $f:M\to X$ is an isometric immersion. If for
each $p\in M$ and each plane $\sigma\subset\nu_pM$ the sectional curvature
$K(\sigma)$ of $X$ vanishes we abbreviate this  by $K_X(\nu M)=0$.

\bigskip\noindent{\bf Lemma 5.3} {\sl Let $\pi:V\to X\ ,\ f:M\to X$ and $\hat f:\hat
M\to V$ be as above, $p\in X$ and $U\subset T_pX$ a linear subspace.

(i) If $K(\sigma)=0$ for all $2$-dimensional subspaces of $U$  then $\exp:U\to X$
 is a totally geodesic
isometric immersion. In particular $K_X(\nu M)=0$ is equivalent to $M$ having
flat sections.

(ii) If $\nu M$ is flat and $K_X(\nu M)=0$ then $\nu\hat M$ is flat as
well. In fact, the horizontal lift of a parallel normal field along an open
subset of $M$ is parallel in $\nu\hat M$.}

\bigskip\noindent{\bf Proof}:

(i) Let $\hat p\in\pi^{-1}(p)$ and $\hat U\subset T_{\hat p}V$ be
the horizontal lift of $U$ through $\hat p$. It suffices to show that the
affine subspace  $\tilde U:=\hat p+\hat U$ meets fibres always orthogonally and
hence projects isometrically into $X$. Let $\hat q\in \tilde U$ and $v$ be a tangent
vector to the fibre at $\hat q$. Then there exists a variation of the geodesic
$c:[0,1]\to V$ from $\hat p$ to $\hat q$ by horizontal geodesics which cover
$\pi\circ c$ and whose variational vector field $Y$ satisfies $Y(1)=v$. The
orthogonal
projection of this Jacobi field $Y$ onto the totally geodesic subspace $\tilde U$ is
a Jacobi field again.
It follows from O'Neill's formula for the curvature of horizontal planes that
the O'Neill tensor vanishes on $\hat U$ and that therefore $<Y'(0),\hat
w>=-<Y(0),\nabla_{\dot c(0)}\hat W>=0$ for all $\hat w\in\hat U$ where $\hat
W$ is a horizontal extension of $\hat w$. Note that $Y(t)$ is vertical. Hence
the projection $\bar Y$ of $Y$ onto $\tilde U$ satisfies $\bar Y(0)=\bar
Y'(0)=0$ and thus vanishes identically. In particular $v=Y(1)$ is orthogonal
to $\tilde U$.

(ii) follows from the vanishing of the O'Neill tensor $O'N$ on $\nu\hat M$ by
O'Neill's curvature formula, cf. chapter 3. \hfill\qed

\bigskip\noindent{\bf Theorem 5.4}
{\sl Let $\pi:V\to X$ be a Riemannian submersion with minimal fibres from
a separable Hilbert space $V$ onto a finite dimensional manifold $X$.
Let $M$ be a properly immersed isoparametric
submanifold of $X$ with flat sections and $\hat M$ its lift to $V$.
 Then $\hat M$ is an immersed
 proper Fredholm isoparametric submanifold of $V$ of finite codimension.}

\bigskip\noindent{\bf Remarks}:

\noindent
(i) Under a mild restriction on $X$, also the converse
is true: $M$ is isoparametric with flat sections if $\hat M$ is isoparametric (Theorem 6.5).
(ii) If $X$ is analytic, the minimality assumption on the
 fibres is not necessary. (Remark 2 following Theorem 6.5).

\bigskip\noindent{\bf Proof}:  By Lemma 5.1, $\hat M$ is proper Fredholm, and by
Lemma 5.3, $\nu(\hat M)$ is flat and parallel normal fields along $\hat M$ are (at least
locally) horizontal lifts of parallel normal fields along $M$. Therefore each
(local)
parallel manifold $\hat M_{\hat\xi}$ of $\hat M$ is locally the lift of a
parallel manifold $M_\xi$ of $M$.
Let $\Phi_t$ and $\hat\Phi_t$ be the goedesic flows of $X$ and $V$,
respectively and $\xi^*=\Phi_1(\xi),\hat\xi^*=\hat\Phi_1(\hat\xi)$ the
corresponding normal fields along $M_\xi$ and $\hat M_{\hat\xi}$. Then $\hat
M_{\hat\xi}$ is regularizable by Lemma 5.2 and $tr\ \hat A_{\hat\xi^*}=tr\
A_\xi$ at corresponding points. Thus $tr\ \hat A_{\hat\xi^*}$ is constant as
well along $\hat M_{\hat\xi}$ and $\hat M$ is isoparametric.\hfill\qed

\bigskip

Isoparametric submanifolds are defined in purely local terms. Therefore one
can not expect that they give rise in general to a global foliation of the
ambient manifold $X$ by parallel manifolds, even if $X$ is the base of a
Riemannian submersion $\pi:V\to X$ with minimal fibres and the submanifold is
embedded. A simple counterexample is furnished by a small geodesic circle on a
flat $2$-torus. Nevertheless there are many interesting situations where the
submanifold defines a global foliation. To study these we make the following definition
which extends the notion of a parallel manifold of an isoparametric submanifold with
 globally flat normal bundle.

\bigskip\noindent {\bf Definition} Let $X$ be a complete Riemannian manifold and
$f:M\to X$ be an immersed isoparametric submanifold.

(i) If $p\in M$ and $\xi\in\nu_pM$ we let $M_\xi:=\{\tilde\xi\in\nu
M\mid\tilde\xi$ is obtained from $\xi$ by parallel translation in $\nu M$
along any curve in $M$ starting at $p\}$ and $f_\xi:M_\xi\to X$ the
restricition of $\exp$ to $M_\xi$. We call $f_\xi(M_\xi)$, or more precisely
$f_\xi:M_\xi\to X$, a parallel manifold of $M$.

(ii) We say that $M$ defines a global foliation of $X$, if $X$ is the disjoint
union of the parallel manifolds of $M$, that is if
$X={\bigcup\limits_\xi}f_\xi(M_\xi)$ and if for all $\xi_1,\xi_2\in\nu M$
either $f_{\xi_1}(M_{\xi_1})=f_{\xi_2}(M_2)$ or $f_{\xi_1}(M_{\xi_1})\cap
f_{\xi_2}(M_{\xi_2})=\phi$. In this case we call the $f_\xi(M_\xi)$ also the
leaves of the isoparametric foliation.

\bigskip\noindent{\bf Remarks}:

\noindent
(i) If $f:M\to X$ defines a global foliation, then $f(M)$ has no
self-intersections as locally a neighborhood of $f(M)$ is covered by parallel
manifolds. This argument shows more precisely that $f$ factorizes over an
injective immersion. Hence $f(M)$ is a closed embedded submanifold if in
addition $f$ is proper. But even then the parallel manifolds need not be
embedded (see below).

\noindent
(ii) If the isoparametric immersion $f:M\to X$ is proper, then
$X=\cup f_\xi(M_\xi)$ is satisfied automatically as each $x\in X$ has a
closest point on $f(M)$ and hence lies on a geodesic starting orthogonally from $M$.

\bigskip\noindent{\bf Theorem 5.5} {\sl Let $X$ be a finite dimensional
Riemannian manifold which is the base of a Riemannian submersion $\pi:V\to X$
with minimal fibres, where $V$ is a separable Hilbert space. Let $M$ be a
connected Riemannian manifold and $f:M\to X$ a properly immersed isoparametric
submanifold with flat sections.

(i) If the induced mapping $f_*:\pi_1 M\to\pi_1 X$ on the fundamental
groups is surjective (in particular if $X$ is simply connected), then $M$
defines a global foliation of $X$.

(ii) If $X$ is simply connected then moreover $f(M)$ has globally flat
normal bundle and all leaves of the foliation are closed embedded
submanifolds.}

\bigskip\noindent{\bf Proof}: We begin with the proof of (ii) and thus assume
that $X$ is simply connected. It follows from the exact homotopy sequence of
the fibration $\pi: V\to X$ that the fibres are connected. Hence also the lift
$\hat M$ of $M$ to $V$ is connected. By Theorem 5.4 and Corollary 4.3 $\hat M$ is
isoparametric in the sense of Terng and thus gives rise to an isoparametric
foliation of $V$ by closed embedded submanifolds. Moreover the normal bundle
of $\tilde M:=\hat f(\hat M)$ is globally flat. Since the horizontal lift of a
normal vector $\xi\in\nu_p M$ to the fibre is parallel  in $\nu\tilde M$ due
to the vanishing of the O'Neill tensor on $\nu\tilde M$ and since the fibres
are connected, the restriction of a parallel normal field along $\tilde M$ to
a fibre is the horizontal lift of a normal vector to $M$. Therefore the
parallel manifolds of $\tilde M$ contain with each point the whole fibre and
the foliation of $V$ can be pushed down to a foliation of $X$ by closed
embedded submanifolds. These are $f(M)$ and its parallel manifolds since a
parallel normal field along $\tilde M$ projects down to a parallel normal
field along $f(M)$, which also shows that $f(M)$ has globally flat normal
bundle.

(i) In the general situation let $\varrho:\tilde X\to X$ be the
universal cover of $X$ and $\Gamma$ its group of deck transformations.
Furthermore let $f^*:M^*\to\tilde X$ be the lift of $f$, where
$M^*=\{(p,\tilde x)\in M\times\tilde X\mid f(p)=\varrho(\tilde x)\}$ and
$f^*(p,\tilde x)=\tilde x$. Then $f^*$ is proper as well and we conclude from
the surjectivity of $\pi_1M\to\pi_1X$ that $M^*$ is also connected. Since also
$\pi:V\to X$ can be lifted to $\tilde X$, the above results imply that $M^*$
defines a foliation of $\tilde X$ by closed embedded submanifolds. Since
$f^*(M^*)$ is invariant under $\Gamma$, each $\gamma\in\Gamma$ permutes the
parallel manifolds of  $M^*$. Hence the projections of these parallel
manifolds to $X$ (which are the parallel manifolds of $M$) are either disjoint
or equal.\hfill\qed

\bigskip\noindent{\bf Remarks}:

\noindent
(i) An example for the situation described in (i) of Theorem 5.5 (with $X$ not
simply connected) is a totally geodesic $\R P^{n-1}\subset\R P^n$.

\noindent
(ii) In general the parallel manifolds of $M$ in the Theorem above are not
embedded submanifolds (although $f(M)$ is always closed and embedded by the
remark before the Theorem). A counterexample is given by $X:=\R^3_{/\Gamma}$,
where $\Gamma$ is the group generated by a glide rotation around the $z$-axis
and $M$ is the image of the $z$-axis in $X$. If the angle of rotation is an
irrational multiple of $2\pi$ then the parallel manifolds of $M$, which are
geodesics, are neither closed nor embedded submanifolds (except $M$ itself).
However, this phenomenon can not occur if $X$ is compact, as we will see below.

\bigskip
Let $f:M\to X$ be a proper isoparametric immersion with flat sections, where
$X$ is the base of a Riemannian submersion $\pi:V\to X$. By part (i) of Lemma
5.3, the exponential maps $\exp:\nu_pM\to X$ are totally geodesic isometric
immersions. Thus the sections are globally defined, but again they are neither
closed nor embedded submanifolds in general, even if $X$ is a compact simply
connected symmetric space. A counterexample is given by the cohomogeneity$-1$ action
 on $S^2\times S^2$ described in section 2. However,  we have the following weak compactness
result for sections which will be sufficient for later applications.

\bigskip\noindent{\bf Theorem 5.6} {\sl Let $X$ be a compact Riemannian
manifold which is the base of a Riemannian submersion $\pi:V\to X$ with
minimal fibres, where $V$ is a separable Hilbert space. Let $M$ be a
connected and properly immersed isoparametric submanifold of $X$ with flat
sections which defines a global foliation of $X$. Then there exists for each $p\in M$
a discrete cocompact group $\Lambda$ of isometries of $\nu_pM$ such that for
any $\xi_1,\xi_2\in\nu_pM\ ,\ \exp\xi_1$ and $\exp\xi_2$ lie on the same leaf iff
$\xi_1$ and $\xi_2$ lie on the same $\Lambda$-orbit. In particular, the
projection of a section onto the space of leaves factorizes over a compact
quotient of the section.}

\bigskip\noindent{\bf Proof}: By the remarks before Theorem 5.5 we may
assume that $M$ is a compact embedded submanifold of $X$. We fix $p\in M$ and
call $\xi_1,\xi_2\in\nu_pM$ equivalent if $\exp\xi_1$ and $\exp\xi_2$ lie on
the same leaf. We denote by $A$ the equivalence class through $0\in\nu_pM$, i.e.
$A=\{\xi\in\nu_pM\mid\exp\xi\in M\}$. Let $\Lambda$ be the group of isometries
of $\nu_pM$ which preserve the equivalence relation. We first prove that
$\Lambda$ acts transitively on $A$. Note that $\exp(\nu_pM)$ meets the
parallel manifolds of $M$ always orthogonally as can be seen by lifting the
foliation to the Hilbert space $V$. Therefore we may identify the tangent
space of $\nu_pM$ at $\xi\in A$ with the normal space of $M$ at $\exp\xi$.
Now, if $\xi_1,\xi_2\in A$ and $c$ is any piecewise differentiable curve in
$M$ from $\exp\xi_1$ to $\exp\xi_2$ we let $\lambda:\nu_pM\to\nu_pM$ be the
affine isometry with $\lambda(\xi_1)=\xi_2$ and $d\lambda_{\xi_1}=\alpha$
where $\alpha$ denotes the parallel translation along $c$ in $\nu M$. By the
very definition of $\lambda$ and the fact that $\exp:\nu_pM\to X$ is a totally
geodesic isometric immersion we get that $\exp(\xi_1+\zeta)$ and
$\exp(\xi_2+\alpha(\zeta))$ lie on the same leaf for any $\zeta\in\nu_pM$.
Thus $\lambda\in\Lambda$ and $\Lambda$ acts transitively on $A$. Moreover the
isotropy group of $\Lambda$ at $0$ contains the normal holonomy group of $M$
at $p$.

The same arguments can be applied to any equivalence class $A'$ containing a
regular element. The regular elements are those normal vectors for which the
corresponding parallel manifold is of highest dimension, and they are dense in
$\nu_pM$ as one can see again by lifting the foliation to the Hilbert space.
Thus $\Lambda$ also acts transitively  on $A'$. If $\xi\in\nu_pM$ there exists
$\zeta\in\nu M$ with $\exp\zeta=\exp\xi$ and $\Vert\zeta\Vert\le d:=\hbox{diam}
X$. The parallel translation in $\nu M$ of $\zeta$ to $p$ along a curve in $M$
also has length at most $d$ and is equivalent to $\xi$. For each regular $\xi$
we therefore find a $\lambda\in\Lambda$ with $\Vert\lambda(\xi)\Vert\le d$
and by continuity the same is true for any $\xi\in\nu_pM$ provided that we
replace $d$ by any $d'>d$. Hence $\Lambda$ acts cocompactly on $\nu_pM$. It
also acts properly  discontinuously as  there exists, by the compactness of
$M$, an $\epsilon>0$ such that $\xi$ is a diffeomorphism on the
$\epsilon$-tube in $\nu M$, and thus the $\epsilon$-balls in $\nu_pM$ around
each $a\in A$ are disjoint. Since the action of $\Lambda$ is proper and the
regular elements are dense, it follows finally that $\Lambda$ acts
transitively on each equivalence class.\hfill\qed

\bigskip\noindent{\bf Corollary 5.7} {\sl Under the assumptions of the
Theorem, the holonomy of $\nu(M)$ is finite and all leaves of the foliation
are compact embedded submanifolds.}

\bigskip\noindent{\bf Proof}: In the course of proof of the Theorem it turned
out that the normal holonomy group of $M$ at $p$ lies in the isotropy group of
$\Lambda$ at $O$. But this is finite as $\Lambda$ acts properly
discontinuously on $\nu_pM$. Hence $M_\xi=\{\tilde\xi\in\nu M\mid\tilde\xi$ is
obtained from $\xi$ by parallel translation in $\nu M$ along curves in $M\}$
is a finite cover of $M$ for each $\xi\in\nu M$. Since we may assume that $M$
is compact it follows that also all leaves are compact and are hence compact
embedded submanifolds by the remark preceding Theorem 5.5.\hfill\qed

\vskip 1 cm

\noindent{\typc 6. Riemannian submersions and focal points}
\nobreak
\bigskip
The results of this chapter will not be used in the sequel. They
rather complement those of the last chapter by proving a converse to Theorem
5.4 and by relating the equifocal submanifolds of Terng-Thorbergsson [TT] to
the isoparametric ones.

Our main tool for this is an elementary but basic lemma about Riemannian
submersions $\pi:\hat X\to X$ (Lemma 6.1). Essentially it says that focal distances
in the direction of a normal vector $\xi$ do not change if one lifts a
submanifold from $X$ to $\hat X$ (and $\xi$ to $\hat\xi)$. This result was known  in
special cases ([O'N2], Theorem 4 and [TT], Lemma 5.12) but it does not seem
to have been noticed before in this gene\-rality. In case $\hat X$ is a euclidean
space
(possibly an infinite dimensional Hilbert space) the lemma immediately relates the
eigenvalues of the shape operators of the lifted manifold to focal
distances of the submanifold in $X$. In particular the minimality of the
fibres of a Riemannian submersion $\pi:V\to X$, where $V$ is a Hilbert space,
is equivalent to a special behaviour of conjugate points in $X$:
For any geodesic $c$ in $X$ the points conjugate to $c(0)$ along $c$ have to
lie symmetrically on both sides of
$c(0)$, at least in the average. This
property is obviously satisfied for symmetric
spaces but also more generally for spaces in which geodesics are orbits of
$1$-parameter groups like normal homogeneous spaces.

Recall that the differential  of the normal exponential map
$\exp_{\nu M}:\nu M\to X$  of a submanifold $M\subset X$ at
$\xi\in\nu_pM$ may be described by means of $M$-Jacobi fields along
$c(t):=\exp t\xi$. These
are by definition the variational vector fields of variations of $c$ through
geodesics starting perpendicularly from $M$. Alternatively, they are the Jacobi
fields $Y$ along $c$ with $Y(0)\in T_pM$ and $Y'(0)^T=-A_\xi
Y(0)$, where $^T$ denotes the tangential component and $A_\xi$ the shape operator
of $M$ in the direction of $\xi$.
If $\J$ denotes the space of $M$-Jacobi fields along $c(t)=\exp t\xi$ then
$T_\xi\nu M$ may be identified with $\J$ by mapping $\dot\xi(0)$ onto
$Y(t):={\partial\alpha\over\partial s}(0,t)$, where $\xi(s)$ is a curve in
$\nu M$ with $\xi(0)=\xi$ and $\alpha(s,t)=\exp t\xi(s)$. Under this
identification,
$(d\exp_{\nu M})_\xi:T_\xi\nu M\to
T_{\exp\xi}$ corresponds to $Y\mapsto Y(1)$. This description
 of $d\exp_{\nu M}$ is also
valid in the infinite dimensional case, i.e. if $M$ is a submanifold of a
Hilbert manifold. A normal vector $\xi$ is called a multiplicity$-m$ focal direction if $\xi$
lies in the domain of $\exp$ and the dimension of the kernel of $d\exp_{\nu
M}$ at $\xi$ is equal to $m$.

\bigskip\noindent{\bf Lemma 6.1} {\sl Let $X$ and $\hat X$ be Riemannian
manifolds of possibly infinite dimension and $\pi:\hat X\to X$ a Riemannian
submersion. Let $M$ be a submanifold of $X$ and $\hat M:=\pi^{-1}(M)$. Then
$\hat\xi\in\nu\hat M$ is a multiplicity$-m$ focal direction of $\hat M$ if
and only if $\pi_*\hat\xi$ is a multiplicity$-m$ focal direction of $M$. In
particular focal distances of $M$ and $\hat M$ in corresponding normal directions are
equal.}

\bigskip\noindent{\bf Proof}: Let $\xi:=\pi_*\hat\xi$ and $\J$ and $\hat \J$ be
the spaces of $M$- and $\hat M$-Jacobi fields along $c(t):=\exp t\xi$ and
$\hat c(t):=\exp t\hat\xi$, respectively. We then get the following
commutative diagram
$$
\tilde\pi_*\matrix{
\hat\J    &{\buildrel\hat f\over\longrightarrow}&T_{\hat q}\hat X\cr
\downarrow&                                     &\downarrow\pi_*\cr
\J         &{\buildrel f\over\longrightarrow}    &T_qX\cr}
$$
where $\hat q=\hat c(1)\ ,\ q=c(1)\ ,\ f$ and $\hat f$ are the evaluation maps
at $t=1$ and $\tilde\pi_*\hat Y(t):=\pi_*(\hat Y(t))$. We want to show that
$\tilde\pi_*$ induces an isomorphism between the kernels of $\hat f$ and $f$,
respectively. But this follows by elementary linear algebra provided we know
that $\tilde\pi_*$ is surjective and $\hat f$ induces an isomorphism between
$\ker\tilde\pi_*=:{\hat\J}^v$ (the space of everywhere vertical $\hat
M$-Jacobi fields) and $\ker\pi_*=:\V_{\hat q}$. The first fact is obvious while
the second comes essentially from the observation that the horizontal lifts of
$c$ define a diffeomorphism between the fibres over $c(0)$ and $c(1)$  (and hence
do not focalize). More precisely, if ${\tilde\J}^v$ denotes the set of Jacobi
fields along $\hat c$ arising from such variations of $\hat c$ by horizontal
lifts of $c$ then ${\tilde\J}^v\subset{\hat\J}^v$. If $Y\in{\hat\J}^v$ with $Y(0)=0$
then also $Y'(0)=0$ as $Y'(0)$ is horizontal and
$<Y'(0),\eta>=-<Y(0),\eta'>=0$ for any horizontal vector field $\eta$ along
$\hat c$. The Jacobi fields in $\hat\J^v$ are therefore completely determined
by their initial value and thus ${\tilde\J}^v={\hat\J}^v$. Since $\hat
f:{\tilde\J}^v\to\V_{\hat q}$ is an isomorphism, the proof of the lemma is
completed.\hfill\qed

\bigskip
As in the last chapter we denote by $V$ always a real separable Hilbert space
and by $\pi:V\to X$ a Riemannian submersion onto a finite dimensional
manifold. If $N\subset V$ is a submanifold, $\lambda\in\R\setminus\{0\}$, and
$\xi\in\nu N$,\quad then $1/\lambda$ is a multiplicity$-m$ focal direction of $N$ iff
$\lambda$ is a multiplicity$-m$ eigenvalue of the shape operator $A_\xi$ of
$N$ (This is simply because $N$-Jacobi fields in $V$ are of the form $A+tB$
with $A,B$ parallel, $A(0)$ tangential and $B(0)^T=-A_\xi A(0)$).

\bigskip
\noindent{\bf Corollary 6.2} {\sl Let $\pi:V\to X$ be a Riemannian submersion,
 $M\subset X$ a submanifold and $\hat M=\pi^{-1}(M)$.
Furthermore let $\hat\xi\in\nu\hat M$ and $\xi:=\pi_*\hat\xi$. Then
$\lambda\in\R\setminus\{0\}$ is a multiplicity$-m$ eigenvalue of the shape
operator $A_{\hat\xi}$ of $\hat M$  iff ${1/\lambda}\cdot\xi$ is a multiplicity$-m$
focal direction of $M$.}

\bigskip
In case $M$ is a point, the corollary relates the eigenvalues of the shape operator of a
fibre to conjugate points in $X$.
In particular one gets the following result.

\bigskip
\noindent{\bf Corollary 6.3} {\sl Let $\pi:V\to X$ be a Riemannian submersion.
If the geodesics of $X$ are orbits of $1$-parameter groups, in particular if
$X$ is a symmetric or more generally a normal homogeneous space, then $\pi$
has minimal fibres.}

\bigskip
\noindent{\bf Proof}: Let $c(t):=\varphi_t(p)$ where $\varphi_t$ is a
$1$-parameter group of isometries. If $c(t)$ is conjugate to $c(0)$ along $c$
then also $c(t+s)=\varphi_s(c(t))$ is conjugate to $c(s)=\varphi_s(c(0))$ and
of the same order. Choosing $s=-t$ shows that $c(-t)$ is conjugate to $c(0)$.
By Corollary 6.2 the eigenvalues of the shape operators of the fibres of $\pi$
therefore come in pairs of opposite sign and thus cancel in the
trace.\hfill\qed

\bigskip
Another consequence of Corollary 6.2 is the following
estimate for the eigenvalues of the shape operator of $\hat M$.

\bigskip
\noindent{\bf Corollary 6.4} {\sl Let $\pi:V\to X$ be a Riemannian submersion
with $X$ compact (more generally with positive injectivity radius), $M\subset
X$ a submanifold and $\hat M=\pi^{-1}(M)$. If $\mu_1\le \mu_2\le\dots\le
0\le \dots\le\lambda_2\le\lambda_1$ are the eigenvalues of the shape operator
$A_{\hat\xi}$ of $\hat M$ then there exists a constant $C$ with
$\lambda_k,\vert\mu_k\vert \le C/k$. In particular
$\Sigma\lambda_k^s+\Sigma\vert\mu_k\vert^s<\infty$ for all $s>1$.}

\bigskip
\noindent{\bf Proof}: Let $\xi:=\pi_*\hat\xi$ and $c(t):=\exp t\xi$. By
Corollary 6.2 it suffices to show that the number of focal points (counted
with multiplicity) of $M$ along $c$ in $[0,T]$ (and $[-T,0])$ grows at most
linearly with $T$. But this number can be estimated by the index of
$c_{\vert_{[0,T]}}$ and this in turn by the dimension of broken Jacobi fields
along $c_{\vert_{[0,T]}}$ where the distance of break points is less then the
injectivity radius of  $X$ (cf. [BC], [M]).\hfill\qed

In generalization of [TT] we call an immersed submanifold $M$ of $X$ {\sl
equifocal} if $\nu M$ is flat and if focal distances and multiplicities in
direction of a normal vector are invariant under parallel translation of
this vector in $\nu M$ along any curve. Actually Terng and Thorbergsson only considered
equifocal manifolds with globally flat and abelian normal bundle in
compact symmetric spaces whose metric is induced by the Killing form.

\bigskip\noindent{\bf Theorem 6.5} {\sl Let $\pi:V\to X$ be a Riemannian
submersion with minimal fibres, $M$ a properly immersed submanifold of $X$ and
$\hat M$ its lift to $V$.  If $X$ admits a compact quotient, then the following
statements are equivalent.

(i) $M$ is isoparametric with flat sections,

(ii) $M$ is equifocal with flat sections, and

(iii) $\hat M$ is isoparametric.

\noindent
The equivalence of (ii) and (iii) does not require
the minimality of the fibres.}

\bigskip
\noindent{\bf Proof}: The implication (i) $\Rightarrow$ (iii) is the content
of Theorem 5.4. If $M$ is equifocal with flat sections, then we conclude from Lemma 5.3 and
Corollary 6.2 - which do not use the minimality of the fibres - that $\nu\hat
M$ is flat and that the non-zero eigenvalues of the shape operator
$A_{\hat\xi}$ of $\hat M$ are constant under parallel translation of
$\hat\xi$. As in the proof of Theorem 4.2 it follows that $\hat M$ is isoparametric.
Thus (ii) implies (iii) as well. If $\hat M$ is isoparametric and $M$ has flat sections
 then one can reverse the above reasoning to get (i) and (ii) from (iii).

Hence we are left to show that $M$
has flat sections if $\hat M$ is isoparametric and $X$ admits a compact
quotient. We thus make these two assumptions and by lifting $\pi:V\to X$ to
the universal cover of $X$ we may assume in addition that $X$ is simply
connected. We also may assume that $M$ and hence $\hat M$ are connected and
that $M$ is embedded. It then suffices to show that the horizontal lift of a
normal vector to $M$ is a parallel normal field along the fibre since then
the O'Neill tensor vanishes on $\nu\hat M$ and thus $K_X(\nu M)=0$.

Let $W$ be the affine span of $\hat M$ in $V$ and $\hat\xi_0\in\nu\hat M$ with
$\hat\xi_0$ tangent to $W$. We denote by $\hat\xi$ the normal field along the
fibre which is the basic horizontal lift of $\pi_*\hat\xi_0$ and by
$\tilde\xi$ the parallel normal field along the fibre extending $\hat\xi_0$.
The eigenvalues of $A_{\hat\xi}$ and $A_{\tilde\xi}$ are both constant because
of Corollary 6.2 and since $\hat M$ is isoparametric. But these eigenvalues are given
by $<\hat\xi,n_i>$ and $<\tilde\xi,n_i>$,
respectively, where the $n_i$ are the curvature normals. Hence
$\hat\xi-\tilde\xi$ is perpendicular to the $n_i$ and therefore also to $W$.
Since $\tilde\xi$ is always tangent to $W$ and
$\hat\xi=\tilde\xi+(\hat\xi-\tilde\xi)$ has the same length as $\tilde\xi\ ,\
\hat\xi$ and $\tilde\xi$ necessarily coincide. The parallel manifold $\hat
M_{\tilde\xi}$ contains therefore in particular with each point the whole fibre through that
point and the same is true for $W$ which is a union of parallel manifolds.

It follows from the Cheeger-Gromoll splitting theorem that $X$ splits as
$X_1\times \R^k$ with compact $X_1$. Since $X_1\times\{x\}$ has no focal
points in $X$ and $X$ is simply connected, the inverse images
$\pi^{-1}(X_1\times\{x\})$ are connected totally geodesic submanifolds and
hence parallel affine subspaces. Thus we may assume that $V=V_1\times\R^k$ and
$\pi((v_1,x))=(\bar\pi v_1,x)$, where $\bar\pi:V_1\to X_1$. The intersection of $W$
(the affine span of $\hat M$) with $V_1\times\{x\}$ is
either empty or $V_1\times\{x\}$. In fact $W\cap(V_1\times\{x\})$ is an affine
subspace which contains with each point the whole fibre of $\pi$ and due to
the compactness of $X_1$ these fibres span $V_1\times\{x\}$ as an affine
subspace. Hence $W$ is of the form $W=V_1\times U$ where $U$ is an affine
subspace of $\R^k$. But this immediately implies that also the parallel
translation of  a $\hat\xi_0\in\nu\hat M$, which is perpendicular to $W$,
coincides with the basic horizontal lift of $\pi_*\hat\xi_0$.\hfill\qed

\bigskip
\noindent{\bf Remarks}:

\noindent
1. Without the compactness assumption on $X$, (iii) does not imply (i)
or (ii) in the Theorem. A counterexample is given by taking $V$ to be $\R^3\
,\ X$ the quotient of $\R^3$ by a $1$-parameter family of glide rotations
around an axis $A$ and $M$ the point in $X$ corresponding to $A$. Then $\hat M$
 is an isoparametric submanifold but $M$ does not have flat
sections.

\noindent
2. In an analytic Riemannian manifold isoparametric submanifolds are
always equifocal. In fact, the mean curvature of the parallel manifolds is
given along $\exp t\xi$ as the logarithmic derivative of the volume distortion
$v(t):=\vert\det d\tilde{\exp}_{t\xi}\vert$, where $\tilde{\exp}$ denotes the
normal exponential map. If this logarithmic derivative is invariant under
parallel translation of $\xi$ for small $t$ then $v(t)$ itself is invariant
under parallel translation of $\xi$ for small $t$ and hence in the analytic
case for all $t$. But the zeros of $v(t)$ are precisely the focal distances.
 Hence (i) implies (iii) without the minimality assumption on the fibres, if
 $X$ is analytic.

\vskip 1 cm

\noindent{\typc 7. A Chevalley-type restriction theorem}
\nobreak
\bigskip
Let $X$ be a compact Riemannian manifold which is the base of a Riemannian
submersion $\pi:V\to X$ with minimal fibres, where $V$ is a separable Hilbert
space, and let $M$ be a properly immersed isoparametric submanifold of $X$ with
flat sections which defines a global foliation of $X$. (A typical example is  $X$ a
compact Lie group
with biinvariant metric and $M$ the conjugacy class of a regular element).
Furthermore let $p\in M$ and $\Sigma:=\exp(\nu_pM)$ a section (which in
general will be only the image of a totally geodesic isometric immersion, see
Lemma 5.3). For any space $A$ of functions on $X$ or $\Sigma$ we denote by
$A^M$ the subspace of functions which are constant on the leaves of the
foliation (the intersection of the leaves with $\Sigma$, respectively). Finally we
denote for any Riemannian manifold $N$
by $C^\Delta(N)$ the subspace of
$C^\infty$-functions which consists finite sums of eigenfunctions of the
Laplacian. For example, $C^\Delta(S^n)$ is the space of
restrictions of polynomials on $\R^{n+1}$ to $S^n$, and is hence a finitely
generated algebra. In case of
the isometrically immersed $\Sigma$ above we interpret $C^\Delta(\Sigma)$ as
$\{f:\Sigma\to\R\mid f\circ\exp_{\mid_{\nu_pM}}\in C^\Delta(\nu_pM)\}$, which
of course coincides with the space of finite sums of eigenfunctions of the
Laplacian on $\Sigma$, if $\Sigma$ is embedded.

\bigskip\noindent{\bf Remark}: If $G$ is a group of isometries of the
Riemannian manifold $N$, then a function $f\in C^\infty(N)$ is called
$G$-finite if the vector space spanned by the functions $f\circ g\ ,\ g\in G$,
is finite dimensional. If $G$ is in addition compact and transitive, then it
is not hard to see that $C^\Delta(N)$ coincides with the space of $G$-finite
functions on $N$. If moreover $N$ is $G$-equivariantly embedded into some
euclidean space (which is always possible by a result of D. Moore) then it has
been shown by Helgason [Hel1] that the space of $G$-finite functions and hence
$C^\Delta(N)$ coincides with the space of restrictions of polynomials from
euclidean space to $N$, just as in the case of the sphere.

\bigskip
A major goal of this section is the proof of the following Chevalley-type
restriction theorem.

\bigskip\noindent{\bf Theorem 7.1} {\sl Under the assumptions above
$C^\Delta(X)^M$ and $C^\Delta(\Sigma)^M$ are isomorphic. In fact, an
isomorphism is given by restricting each $f\in C^\Delta(X)^M$ to $\Sigma$.}

\bigskip
Note that we may identify $C^\Delta(\Sigma)^M$ with
$C^\Delta(\nu_pM)^\Lambda$ where $\Lambda$ is the discrete cocompact group of
isometries of $\nu_pM$ introduced in Theorem 5.6. Hence $C^\Delta(\Sigma)^M$
may be viewed as the space of finite sums of eigenfunctions on a flat torus
which are invariant under a finite group of isometries, namely the quotient
of $\Lambda$ by a lattice.

To prepare the proof of the theorem we begin with several lemmas.
Let $\Gamma\subset\R^n$ be a lattice. Then $C^\Delta(\R^n)^\Gamma
=\{{\sum\limits_{\omega\in\Gamma^*}}a_\omega e^{2\pi i\omega}\mid a_\omega\in
\C\ ,\ a_{-\omega}=\bar a_\omega\ , \ a_\omega=0$ except for
finitely many $\omega\}$, where $\Gamma^*:=\{\omega:\R^n\to\R\mid\omega$
linear, $\omega(\Gamma)\subset\Z\}$ denotes the dual lattice. If
$f={\sum\limits_{\omega\in\Gamma^*}}a_\omega e^{2\pi i\omega}$, we call
$\hbox{supp}_\Gamma(f):=\{\omega\in\Gamma^*\mid a_\omega\ne 0\}$ the $\Gamma$-support
of $f$. The following lemma will be crucial.

\bigskip\noindent{\bf Lemma 7.2} {\sl Let $\alpha\in\Gamma^*$ and $f\in
C^\Delta(\R^n)^\Gamma$. If $f\tan\pi\alpha$ or $f\cot\pi\alpha$
extends to a $C^\infty$-function $g$ on $\R^n$, then $g\in
C^\Delta(\R^n)^\Gamma$ and $\hbox{supp}_\Gamma(g)$ is contained in the convex hull of
$\hbox{supp}_\Gamma(f)$.}

\bigskip\noindent{\bf Proof}: We only consider $g=f\cot\pi\alpha$, the other
case being similar. Let $f=\Sigma a_\omega e^{2\pi i\omega}$. Since $g$ is
$\Gamma$-invariant, we may expand it also into a Fourier series
$g={\sum\limits_{\omega \in\Gamma^*}}b_\omega e^{2\pi i\omega}$. From
$$
(e^{2\pi i\alpha}+1)\Sigma a_\omega e^{2\pi i\omega}=(e^{2\pi i
\alpha}-1)\Sigma b_\omega e^{2\pi i\omega}
$$
we get $a_{\omega-\alpha}+a_\omega=b_{\omega-\alpha}-b_\omega$ and hence
$$
b_{\omega-\alpha}=b_\omega+a_{\omega-\alpha}+a_\omega \eqno(*)
$$
or equivalently
$$
b_{\omega+\alpha}=b_\omega-(a_\omega+a_{\omega+\alpha})\eqno(**)
$$
for all $\omega\in\Gamma^*$.

If $\omega_0$ is not contained in the convex hull of $\hbox{supp}_\Gamma(f)$ then this is also
the case either for all $\omega_0+k\alpha$  or for all
$\omega_0-k\alpha$ with $k\in\N_0$. In the first case we conclude from $(**)\quad
b_{\omega_0}=b_{\omega_0+\alpha}=b_{\omega_0+2\alpha}=\dots$ and hence
$b_{\omega_0}=0$, since $g$ lies in $L^2$ when considered as a function on $\R^n/\Gamma$. Similarly
we conclude in the second case from $(*)\quad
b_{\omega_0}=0$, which completes the proof.\hfill\qed

\bigskip\noindent{\bf Lemma 7.3} {\sl Let $\Lambda\subset I(\R^k)$ be a discrete
cocompact group of isometries of\quad $\R^k$ and $W\subset \Lambda$ a subgroup generated
by reflections in a  $\Lambda$-invariant family $\H$ of affine hyperplanes.
Then $W$ does not split off a finite factor and is thus
 an affine Weyl group (whose action on $\R^k$, however,
might contain a trivial summand). In particular each $H\in\H$ is contained in an infinite family
of parallel affine hyperplanes of $\H$.}

\bigskip\noindent{\bf Proof}: From the theory of affine reflection groups we
have (after a possible adjustment of the origin) a decomposition of
$\R^k$ into $V_0\oplus V_1\oplus V_2$ and of $W$ into $W_1\times W_2$ where $W_1$ is an affine
Weylgroup acting only on $V_1,W_2$ is a finite linear reflection group acting
only on $V_2$ and $V_0$ is the set of fixed vectors of $W$. The reflection
hyperplanes of $W_2$ are precisely those which have no further parallel
hyperplanes in $\H$. This subset of hyperplanes is left invariant by $\Lambda$ and
hence also their intersection, which is equal to $V_0\oplus V_1$. Since $\Lambda$ acts
cocompactly, $V_0\oplus V_1=\R^k$ and thus $W_2=\{e\}$.\hfill\qed

\bigskip
\noindent{\bf Lemma 7.4} {\sl Let $\Lambda\subset I (\R^k)$ be a discrete
cocompact group of isometries of $\R^k$. Then $C^\Delta(\R^k)^\Lambda$ is dense in
$C^\infty(\R^k)^\Lambda$ with respect to the supremum norm.}

\bigskip\noindent{\bf Proof}: Let $\Gamma\subset \Lambda$ be the subgroup of
translations in $\Lambda$. Then $\Gamma$ is a lattice and $\bar \Lambda:=\Lambda/\Gamma$ is a
finite group which operates on $\R^k/\Gamma$. Since $\R^k/\Gamma$ is compact,
$C^\Delta(\R^k/\Gamma)$ is dense in $C^\infty(\R^k/\Gamma)$ with respect to
the supremum norm. Hence also $C^\Delta(\R^k/\Gamma)^{\bar \Lambda}$ is dense in
$C^\infty(\R^k/\Gamma)^{\bar \Lambda}$ as follows by averaging an approximating
sequence in $C^\Delta(\R^k/\Gamma)$ over $\bar \Lambda$. From this the lemma follows immediately.\hfill\qed

Our final preparation for the proof of the theorem is an extension of remark
5.1 of [KT] in which the mean curvature vector field of an isoparametric {\sl
hypersurface} in a Hilbert space had been given.

Let $\hat M\subset V$ be a connected, proper Fredholm isoparametric
submanifold of finite codimension of a separable Hilbert space $V\ ,\ p\in\hat
M$ and $\hat\Sigma:=p+\nu_p\hat M$ be the affine normal space at $p$. Let $\H$
denote the set of focal hyperplanes of $\hat M$ in
$\hat\Sigma,\hat\Sigma_{reg}:=\hat\Sigma\setminus{\bigcup\limits_{H\in\H}}H$ the set of
regular points and $W$ the group generated by reflections in the hyperplanes
of $\H$. To each $H\in\H$ there is associated its so called multiplicity $m_H$
which is equal to the dimension of the curvature sphere of $\hat M$ through
$p$ focalizing on $H$. This number is invariant under $W$, that is
$m_{w(H)}=m_H$ for all $w\in W$, and it is also independent of
$p\in\hat\Sigma_{reg}$. More precisely, for any $q\in\hat\Sigma_{reg}$ the
affine normal space of the parallel manifold $\hat M_\xi$ of $\hat M$ through $q$ is equal to
$\hat \Sigma$ and the focal hyperplanes of $\hat M_\xi$ in that affine normal
space and their
multiplicities coincide with that of $M$. It is well known that $\H$ contains finitely many hyperplanes,
say $H_1,\dots,H_r$, such that each
 $H\in \H$ is parallel to  one of these and the $H_i$ have a point in common.
We endow $\hat\Sigma$
with the structure of a euclidean vector space by choosing as origin a point
from the intersection of the $H_i$. Let $v_i\in\hat\Sigma$ be a unit vector
perpendicular to $H_i$. We assume that $\H$ contains with each
$H\in\H$ another parallel hyperplane and let $l_i>0$ be the common distance
between neighboring parallel hyperplanes in the infinite family determined by $H_i$.
Since the reflection in a hyperplane interchanges its parallel neighbors there
are numbers $m_i^+,m_i^-,i\in\{1,\dots,r\}$, such that $m_H=m_i^+$ if $H=2k\
l_i v_i+H_i$ and $m_H=m_i^-$ if $H=(2k+1)l_iv_i+H_i$ for all $k\in\Z$.

\bigskip\noindent
{\bf Lemma 7.5} {\sl The parallel manifolds of $\hat M$ through regular points
of $\hat\Sigma$ are regularizable. More precisely, the shape operator of the
parallel manifold of $\hat M$
through $q\in\hat\Sigma_{reg}$ in direction
$\xi\in\hat\Sigma$ is equal to $<\hat\eta(q),\xi>$ where $\hat\eta$ is given by
$$
\hat\eta={\sum\limits^r_{i=1}}{\pi\over 2l_i}(m_i^-\tan \pi
\alpha_i-m_i^+\cot \pi \alpha_i)v_i
$$
and $\alpha_i$ is the linear functional on $\hat\Sigma$ with
$\alpha_i(x)=<x,v_i>/2l_i$. Note that $2\alpha_i(x)\in\Z$ iff $x$ lies on a
focal hyperplane parallel to $H_i$, which is in particular the case  if
$x+\H\subset\H$.}

\bigskip\noindent{\bf Proof}: The non-zero eigenvalues of the shape operator
in the direction of $\xi$ are
the numbers $1/t$ with $q+t\xi\in{\bigcup\limits_{H\in\H}}H$ and their
multiplicity is the sum over all $m_H$ for which $q+t\xi\in H$. These
eigenvalues fall naturally into $2r$ classes, namely into the classes of
eigenvalues $1/t$ with multiplicity $m^+_i$ where $q+t\xi\in 2k\ l_iv_i+H_i$,
i.e. $<q+t\xi,v_i>=2k\ l_i$ or equivalently $1/t={<\xi,v_i>\over
2kl_i-<q,v_i>}$ and the classes of eigenvalues $1/t$ with multiplicity $m_i^-$
where $q+t\xi\in(2k+1)l_iv_i+H_i$, i.e. ${1\over
t}={<\xi,v_i>\over(2k+1)l_i-<q,v_i>}\ ,\ k\in\Z$.

The regularized trace for these subclasses exists, since
${\sum\limits_{k\in\Z}}{1\over z+k}:={1\over
z}+{\sum\limits^\infty_{k=1}}({1\over z+k}+{1\over z-k})=\pi\cot\pi z$ by the
well known partial fraction formula for the cotangent. More precisely
$m_i^+\cdot{\sum\limits_{k\in\Z}}{<\xi,v_i>\over 2kl_i-<q,v_i>}=-{m_i^+<\xi,v_i>\over 2 l_i}
\pi\cot\pi\alpha_i(q)$ and $m_i^-\cdot
{\sum\limits_{k\in\Z}}{<\xi,v_i>\over(2k+1)l_i-<q,v_i>}=m_i^-{<\xi,v_i>\over
2l_i}\pi\tan\pi\alpha_i(q)$, where $\alpha_i(q)={<q,v_i>\over 2l_i}$. Thus the
regularized trace of the shape operator exists by Lemma 4.4 and is equal to the
stated formula.\hfill\qed

\bigskip\noindent{\bf Proof of  Theorem 7.1}: We may assume that $M$ is an
embedded submanifold of $X$. Let $\hat M:=\pi^{-1}M$ and $\hat M_0$ be a
connected component of $\hat M$. Then $\hat M_0$ is a proper Fredholm
isoparametric submanifold of $V$ of finite codimension and $\hat M$ consists of
certain parallel manifolds of $\hat M_0$ which are also regular, i.e. have the
same codimension as $\hat M_0$. Let $\hat p\in\hat M_0$ and $\hat \Sigma:=\hat
p+\nu_{\hat p}\hat M_0$ be the affine normal space.
 According to Terng [T2], the restriction mapping defines an isomorphism
from $C^\infty(V)^{\hat M_0}$ onto $C^\infty(\hat\Sigma)^{\hat M_0}$ (the set
of $C^\infty$-function on $\hat\Sigma$ which are constant on the intersections
of $\hat \Sigma$ with the parallel manifolds of $\hat M_0$) and
$C^\infty(\hat\Sigma)^{\hat M_0}=C^\infty(\hat\Sigma)^W$ where $W$ is the
group generated by reflections in the focal hyperplanes of $\hat M_0$ in
$\hat\Sigma$. This isomorphism yields an isomorphism $C^\infty(V)^{\hat
M}\to C^\infty(\hat\Sigma)^{\hat M}$ where $C^\infty(V)^{\hat M}$ is the set
of $C^\infty$-functions on $V$, which are constant on the preimages under $\pi$
of the parallel manifolds of $M$ in $X$ and $C^\infty(\hat\Sigma)^{\hat M}$
is the set of their restrictions to $\hat\Sigma$. Since we may identify
$C^\infty(X)^M$ with $C^\infty(V)^{\hat M}$  and
$C^\infty(\Sigma)^M$ with $C^\infty(\hat\Sigma)^{\hat M}$, by mapping $f$ to
$f\circ\pi$, the restriction mapping
$C^\infty(X)^M\to C^\infty(\Sigma)^M$ is an
isomorphism, where $C^\infty(\Sigma):=\{f:\Sigma\to\R\mid f\circ\exp\in
C^\infty(\nu_pM)\}=\{f:\Sigma\to\R\mid f\circ\pi\in C^\infty(\hat\Sigma)\}$.
Note that we may identify $\hat\Sigma$ with $\nu_pM$ by means of the isometry
$\hat p+\hat\xi\mapsto\pi_*(\hat\xi)$ and that $\pi:\hat\Sigma\to X$
corresponds to $\exp:\nu_pM\to X$ under this identification. Let $\Lambda$ be
the group of isometries of $\hat \Sigma$ which map each $x\in\hat\Sigma$ to
an equivalent point, i.e. to a point $y$ such that $\pi x$ and $\pi y$ lie on
the same parallel manifold of $M$. Then $C^\infty(\hat\Sigma)^{\hat
M}=C^\infty(\hat\Sigma)^\Lambda$ by Theorem 5.6 and the action of $\Lambda$ on
$\hat\Sigma$ is discrete and cocompact. Moreover $\Lambda$ leaves the set $\H$
of focal hyperplanes of $\hat M_0$ invariant and contains $W$. Hence $W$ is an
affine Weyl group by Lemma 7.3 and each $H\in\H$ belongs to an infinite family
of parallel focal planes.

Let $\delta:C^\infty(\Sigma)^M\to C^\infty(\Sigma)^M$ be defined by
$\delta(f):=(\Delta^X\tilde f)_{\vert_{\Sigma}}$ where $\Delta^X$ denotes the
Laplacian of $X$ and $\tilde f\in C^\infty(X)^M$
 the extension of $f$ to $X$
along the parallel manifolds of $M$. Note that $\Delta^X$
 leaves $C^\infty(X)^M$ invariant. In
fact, the projection from $X$ onto the space of leaves of the isoparametric
foliation is a Riemannian submersion outside the singular leaves,
and the mean curvature vector field of
the leaves is a basic horizontal vector field. If $f\in C^\infty(X)^M$,
$\Delta^Xf$ is therefore constant along the regular leaves, by the standard
formula for the Laplacians of a Riemannian submersion and hence
constant on all leaves by continuity.

Let $N_1,N_2\subset X$ be two submanifolds of $X$ of complementary dimension
which meet orthogonally at $q$. Then
$$
\Delta^Xf(q)=\Delta^{N_1}f_{\vert_{N_1}}(q)+\Delta^{N_2}f_{\vert_{N_2}}(q)+\eta_1(f)+\eta_2(f)
$$
where $\eta_1$ and $\eta_2$ denote the mean curvature vectors of $N_1$ and
$N_2$ at $q$. This formula follows easily from $\Delta^Xf(q)=-tr\
Hess_qf=-{\sum\limits_i}<\nabla_{e_i}grad f,e_i>$ by using an orthonormal
basis $\{e_i\}$ of $T_qX$ compatible with the decomposition $T_qN_1\oplus
T_qN_2$. Specializing it to $N_1=\Sigma$ (which we may assume to be embedded
for this consideration) and $N_2$ a parallel manifold $M_\xi$ of $M$ of
maximal dimension yields
$$
\delta(f)(q)=\Delta^\Sigma f(q)+\eta(f)
$$
for all $f\in C^\infty(\Sigma)^M$, where $q\in\Sigma\cap M_\xi$ and $\eta$ denotes
 the mean curvature of $M_\xi$ at $q$.

We now identify $C^\infty(\Sigma)^M$ with $C^\infty(\hat\Sigma)^\Lambda$ by
means of $f\mapsto f\circ\pi_{\vert_{\hat\Sigma}}$ and denote the endomorphism
of $C^\infty(\hat\Sigma)^\Lambda$ corresponding to $\delta$ by $\hat\delta$.
Since the regularized mean curvature $\hat\eta$ of the isoparametric
submanifold $\pi^{-1}(M_\xi)$ at $\hat q\in\hat\Sigma\cap\pi^{-1}(q)$
corresponds under $\pi_*$ to $\eta$ by Lemma 5.2, we get
$$
\hat\delta f(\hat q)=\Delta^{\hat\Sigma}f(\hat q)+\hat\eta(f)
$$
for all $f\in C^\infty(\hat\Sigma)^\Lambda$. But for $\hat\eta$ we have the
explicit formula of Lemma 7.5 and thus obtain (using the notation introduced
before Lemma 7.5)
$$
\hat\delta (f)=\Delta^{\hat\Sigma}(f)+{\sum\limits^r_{j=1}}{\pi\over 2
l_j}(m_j^-\tan\pi\alpha_j-m_j^+\cot\pi\alpha_j)v_j(f)
$$
for all $f\in C^\infty(\hat\Sigma)^\Lambda$. Recall that $\hat\Sigma$ is
viewed as a vector space by introducing as origin a point through which passes
a focal hyperplane from each infinite family of parallel focal hyperplanes and
that the $\alpha_j$ are linear functionals on $\hat\Sigma$.

Our aim is to prove that the isomorphism $C^\infty(X)^M\to
C^\infty(\hat\Sigma)^\Lambda\ ,\ f\mapsto f\circ\pi_{\vert_{\hat\Sigma}}$,
maps $C^\Delta(X)^M$ onto $C^\Delta(\hat\Sigma)^\Lambda$. Denoting the inverse
mapping by $\psi$, the essential step will be to show that
$\psi(C^\Delta(\hat\Sigma)^\Lambda)\subset C^\Delta(X)^M$ or equivalently that
each $f\in C^\Delta(\hat\Sigma)^\Lambda$ is a finite sum of eigenfunctions of
$\hat\delta=\psi^{-1}\circ\Delta^X\circ\psi$. This in turn will follow, if each
such $f$ lies in a finite dimensional $\hat\delta$-invariant subspace on which
$\hat\delta$ can be diagonalized. But  the last condition is automatic as
$\hat\delta$ is symmetric on $C^\infty(\hat\Sigma)^\Lambda$ with respect to
the inner product $<f,g>:=<\psi(f),\psi(g)>_{L^2(X)}$.

Let $\Gamma$ be the lattice of translations of $\Lambda$ and
$\Gamma_2:=2\Gamma$. Then $\alpha_j\in\Gamma_2^*$ by Lemma 7.5. Hence
$v_j(f)\tan\pi\alpha_j$ as well as $v_j(f)\cot\pi\alpha_j$ are
$\Gamma_2$-invariant functions. They are smooth in the whole of $\hat\Sigma$
as $\tan\pi\alpha_j$ and $\cot\pi\alpha_j$ only have poles of first order
which moreover are contained in the focal hyperplanes $kl_jv_j+H_j\ ,\
k\in\Z$, and $v_j(f)$ vanishes there due to the invariance of $f$ under
reflections along these hyperplanes. Each $f\in C^\Delta(\hat\Sigma)^\Lambda$
is $\Gamma_2$-invariant and can thus be developed into a Fourier series
${\sum\limits_{\omega\in\Gamma_2^*}}a_\omega e^{2\pi i\omega}$ which is finite
by assumption. Since $e^{2\pi i\omega}$ is an eigenfunction of
$\Delta^{\hat\Sigma}$ and $v_j$, we conclude from Lemma 7.2 and
the above formula for $\hat\delta$ that $\hat\delta$ leaves
$C^\Delta(\hat\Sigma)^\Lambda$ invariant and that moreover the
$\Gamma_2$-support of $\hat\delta f$ is contained in the convex hull of that
of $f$. This shows in particular that $\hat\delta$ leaves invariant the finite
dimensional subspaces $\{f\in C^\Delta(\hat\Sigma)^\Lambda\mid
supp_{\Gamma_2}(f)\subset B_r(0)\}\ ,\ r>0$, which fill out completely
$C^\Delta(\hat\Sigma)^\Lambda$. By the remarks above this implies
$\psi(C^\Delta(\hat\Sigma)^\Lambda)\subset C^\Delta(X)^M$.

Equality follows now by a standard density argument. From Lemma 7.4 we get
that $C^\Delta(\hat\Sigma)^\Lambda$ is dense in $C^\infty(\hat\Sigma)^\Lambda$
with respect to the supremum norm and hence also that
$\psi(C^\Delta(\hat\Sigma)^\Lambda)$ is dense in $C^\infty(X)^M$ and thus
in $C^\Delta(X)^M$. Since $C^\Delta(X)$ is the algebraic direct sum
of the eigenspaces $E_k,k\in\N$, of $\Delta^X$ and since $\Delta^X$ leaves
$C^\Delta(X)^M$ invariant, we have $C^\Delta(X)^M={\bigoplus\limits_k}E^M_k$
(algebraic direct sum) with $E_k^M=E_k\cap C^\infty(X)^M$. From
$\Delta^X\circ\psi=\psi \circ\hat\delta$ it follows that $\Delta^X$ leaves
$\psi(C^\Delta(\hat\Sigma)^\Lambda)$ invariant too, and hence that
$\psi(C^\Delta(\hat\Sigma)^\Lambda)={\bigoplus\limits_k}\bar E_k^M$ for
certain $\bar E_k^M\subset E_k^M$. If $\bar E_k^M\ne E_k^M$ for some $k$,
there would exist an $f\in C^\Delta(X)^M\setminus\{0\}$ perpendicular to
$\psi(C^\Delta(\hat\Sigma)^\Lambda)$ with respect to the $L^2$-inner product.
But this would contradict the density of $\psi(C^\Delta(\hat\Sigma)^\Lambda)$.
Hence $\psi(C^\Delta(\hat\Sigma)^\Lambda)=C^\Delta(X)^M$ and the theorem is
proved.\hfill\qed

\bigskip
\noindent{\bf Theorem 7.6} {\sl Let $X$ be a compact simply connected Riemannian
manifold which is the base of a Riemannian submersion $\pi:V\to X$, where $V$
is a separable Hilbert space and let $M$ be a connected properly immersed
isoparametric submanifold of $X$ with flat sections. Then $M$  defines a
global foliation of $X$ and $C^\Delta(X)^M$ is isomorphic to the
polynomial algebra $\R[X_1,\dots,X_k]$, where $k=\hbox{codim} M$. As free generators
of the algebra $C^\Delta(X)^M$ one may take eigenfunctions of the Laplacian of $X$.}

\bigskip\noindent{\bf Proof}: We retain the notations of the last theorem and
its proof. According to Theorem 5.5, $M$ defines a global foliation of $X$.
In particular we may assume that $M$ is embedded.
Since $X$ is simply connected, $\hat M$ is connected as well and points
$x,y\in\hat\Sigma=\hat p+\nu_{\hat p}\hat M$
are equivalent if and only if they lie on the same parallel manifold of
$\hat M$ and this in turn
is precisely the case if they lie on the same $W$-orbit ([T2], Theorem 9.6).
Therefore $\Lambda=W\ ,\ W$ is an affine Weyl group which acts
cocompactly on $\hat\Sigma$, and $C^\Delta(\Sigma)^M=C^\Delta(\hat\Sigma)^W$. Let $\Gamma$ be the lattice of translations of $W$
and $\bar W$ the isotropy group of $W$ at zero.
Recall that zero is assumed to be a
special point, i.e. a point through which passes
a focal hyperplane from each
family of parallel ones. Then $W$ is the semidirect product of $\Gamma$ with
$\bar W$ and
$C^\Delta(\hat\Sigma)^W\otimes\C=\{{\sum\limits_{\omega\in\Gamma^*}}a_\omega e^{2\pi
i\omega}\mid a_\omega\in\C,a_\omega=0$ except for finitely many $\omega\}^{\bar W}$.
According to [Bou1],
Ch VI, \S 3 Th\'eor\`eme 1, this ring is an algebra, isomorphic to
$\C[X_1,\dots,X_k]$. Moreover the generators of this algebra may be chosen in
a special way as follows.

Since $W$ is an affine Weyl group, there exists a root system $R$ in
the dual of $\hat \Sigma$ such that the reflection hyperplanes are the sets
$H_{\alpha,l}:=\{x\in\hat\Sigma\mid\alpha(x)=l\}$ where $\alpha\in R$ and
$l\in\Z$ (Necessarily $R=\{\pm 2\alpha_i\mid i=1,\dots,r\})$. Let $\Gamma$ as
before denote the lattice of translations of $W=\Lambda$. Identifying
translations $x\mapsto x+v$ of $\hat\Sigma$ with $v\in\hat\Sigma\ ,\Gamma$
coincides with the lattice $Q(\check R)=\{\Sigma m_\alpha\check\alpha\mid
m_\alpha\in\Z,\alpha\in R\}$ of inverse roots (cf. [Bou1], [BD]). Hence
$\Gamma^*$ coincides with the dual lattice $Q(\check R)^*$ of $Q(\check R)$,
which is by definition the lattice of weights.
By choosing a Weyl chamber one
distinguishes a basis $\bar \omega_1,\dots,\bar \omega_k$ of the lattice
$\Gamma^*$  of weights whose elements
are called the fundamental weights. There is a partial order on the set of
weights. In particular $\omega'<\omega$ iff $\omega'$ is contained in the convex hull
of $W\omega$ but not in $W\omega$ itself. Note that $W$ acts on $\Gamma$ and hence on
$\Gamma^*$ as well. If $\omega\in\Gamma^*$ then the so called symmetric sum
$S(e^{2\pi i\omega}):={\sum\limits_{\tilde\omega\in W\omega}}e^{2\pi
i\tilde\omega}$ is an element of $C^\Delta(\hat\Sigma)^W\otimes\C$. Now,
the above quoted theorem states moreover that any elements of the form $S(e^{2\pi
i\bar\omega_j})+
{\sum\limits_{{\omega\in\Gamma^*}\atop{\omega<\bar\omega_j}}}
a_\omega S(e^{2\pi
i\omega}),j=1,\dots,k$, are a basis of the polynomial algebra
$C^\Delta(\hat\Sigma)^W\otimes\C$.

We denote the  $\C$-linear extension of $\hat\delta:C^\Delta(\hat\Sigma)^W\to
C^\Delta(\hat\Sigma)^W$ to the complexification also by $\hat\delta$. This leaves
for all $j\in\{1,\dots,k\}$ the finite dimensional subspaces invariant which
are spanned
by $S(e^{2\pi i\bar\omega_j})$ and the $S(e^{2\pi i\omega})$ with
$\omega<\bar\omega_j$, as $supp_\Gamma\hat\delta f$ is contained in the convex
hull of $supp_\Gamma f$ for all $f\in C^\Delta(\hat\Sigma)^W$.
 Since $\hat\delta$ can be diagonalized on these
subspaces,
there exist eigenfunctions $\varphi_j$ of $\hat\delta$ of the form
$\varphi_j=S(e^{2\pi
i\bar\omega_j})+{\sum\limits_{\omega<\bar\omega_j}}\alpha_\omega S(e^{2\pi
i\omega})$ whose extensions $\tilde\varphi_j:=\psi(\varphi_j)$ to $X$
are thus complex valued
eigenfunctions of $\Delta^X$. If $-\bar\omega_j\in W\bar\omega_j$, then
$S(e^{2\pi i\bar\omega_j})$ is real valued. Eventually replacing
$\varphi_j$ by its real part, which is of the same form $S(e^{2\pi
i\bar\omega_j})+$ lower order terms, we may
assume that $\varphi_j$ and hence $\tilde\varphi_j$ are real valued in this case. In
general however, $-\bar\omega_j$ will not be contained in $W\bar\omega_j$. But
there exists a permutation $\pi$ of $\{1,\dots,k\}$ of order two such that
$-\bar\omega_j\in W\bar\omega_{\pi(j)}$ for all $j=1,\dots,k$ (cf.
[Bou1],[BD]). We may arrange the fundamental weights in such a way that
$\pi(i)=i$ for $i=1,\dots,s$ and $\pi(s+i)=s+i+t$ for $i=1,\dots, t$, where
$s,t$ are nonnegative integers
with $s+2t=k$. Then we may take
$\varphi_1\dots,\varphi_{s+t},\bar\varphi_{s+1},\dots,\bar\varphi_{s+t}$
as a basis for $C^\Delta(\hat\Sigma)^W\otimes\C$ and hence also
$\psi_1,\dots,\psi_k$ where $\psi_i$ denotes the real part of $\varphi_i$ for
$i=1,\dots,s+t$ and
$\psi_{s+t+i}$ denotes the imaginary part of $\varphi_{s+i}$
for $i=1,\dots,t$. Thus $\psi_1,\dots,\psi_k$
are real valued eigenfunctions of $\hat\delta$. Their extensions
$\tilde\psi_1,\dots,\tilde\psi_k$ are therefore real valued
eigenfunctions of $\Delta^X$ and form a
basis of the polynomial algebra $C^\Delta(X)^M$.\hfill\qed

\bigskip\noindent{\bf Remark}: In the non-simply connected case
$C^\Delta(X)^M$ is not a polynomial algebra in general. For example if
$X=\R^n/\Gamma$ is a flat torus and $M$ a point, then
$C^\Delta(X)^M\otimes\C=C^\Delta(X)\otimes\C$ is the space of finite Fourier
series with complex coefficients and thus isomorphic to
$\C[X_1,\dots,X_n,X^{-1}_1,\dots,X_n^{-1}]$.

\bigskip\noindent{\bf Corollary 7.7} {\sl Let $X$ and $M$ be as in the last
theorem. Then there exists a mapping $\psi:X\to\R^k$ whose components
$\psi_i:X\to\R$ are eigenfunctions of the Laplacian of $X$ and whose level sets are
the leaves of the isoparametric foliation determined by $M$.
Moreover $\psi$ is a submersion when
restricted to the set of regular leaves, i.e. to those of maximal dimension.}

\bigskip\noindent{\bf Proof}: We let $\psi:=(\psi_1,\dots,\psi_k)$ with
$\psi_i$ as in the proof of the theorem. Then the $\psi_i$ are eigenfunctions
of $\Delta^X$ and are constant on the parallel manifolds $M_\xi$ of $M$. Since
$C^\infty(\hat\Sigma)^W$ separates the $W$-orbits,
$C^\infty(X)^M$ separates the parallel manifolds of $M$. That is, for any
$M_{\xi_1}\ne M_{\xi_2}$ there exists an $f\in C^\infty(X)^M$ with
$f(M_{\xi_1})\ne f(M_{\xi_2})$.  Due to the density of $C^\Delta(X)^M$ in $C^\infty(X)^M$
with respect to the supremum norm we may choose such an $f$ from
$C^\Delta(X)^M$. As $C^\Delta(X)^M$ is generated by the $\psi_i$ as an algebra, we
find even an index $i$ with $\psi_i(M_{\xi_1})\ne\psi_i(M_{\xi_2})$. Hence
$\psi(M_{\xi_1}) \ne\psi(M_{\xi_2})$ and the level sets of $\psi$ are the
parallel manifolds of $M$.

If $x\in\hat\Sigma_{reg}$ and $v\in T_x\hat\Sigma,v\ne 0$, then there exists
an $f\in C^\infty(\hat\Sigma)^W$ with ${<grad f_x,v>\ne 0}$. In fact, one can
extend any $C^\infty$-function with support in a sufficiently small
neighborhood of $x$ to a $W$-invariant $C^\infty$-function. We may approximate
$f$ and its first derivatives uniformly by functions from
$C^\Delta(\hat\Sigma)^W$. Thus there exists even an $f\in C^\Delta(\Sigma)^W$
with $<grad f_x,v>\ne 0$. The extension $\tilde f$ of $f$ to $X$ lies in
$C^\Delta(X)^M$ and satisfies $\tilde f\circ\pi_{\vert_{\hat\Sigma}}=f$ and hence $<grad\tilde
f_{\pi(x)},\pi_*v>=<grad f_x,v>\ne 0$. As the $\psi_i$ generate
$C^\Delta(X)^M$ as an algebra there exists therefore an index $i$ with
$<(grad\psi_i)_{\pi(x)},\pi_*v>\ne 0$ which implies the surjectivity of
$\psi_*$ at $\pi(x)$.\hfill\qed

\bigskip
In the example $X=G$, a compact Lie group with biinvariant metric, and $M$ a
conjugacy class of a regular element the theorems and the corollary are of
course well known and are weak versions of fundamental results in the theory
of representations of compact Lie groups. In fact, the action of $G$ on itself
by conjugation is hyperpolar with the maximal tori as flat sections. If $T$ is
a maximal torus, then $C^\infty(G)^M=C^\infty(G)^G$ is the set of smooth class
functions on $G$ and $C^\infty(T)^M=C^\infty(T)^{W_G}$, where $W_G$ denotes
the Weyl group of $G$. Moreover, $C^\Delta(G)$ is spanned by the coefficients
of representations of $G$ (Peter-Weyl) and $C^\Delta(G)^G\otimes \C$ is the
complex vector space spanned by the characters of representations of $G$. Thus
$C^\Delta(G)^G\otimes\C$ is isomorphic to $R(G)\otimes \C$, where $R(G)$
denotes the representation ring of $G$. Hence Theorem 7.1 shows in this case
that restriction from $G$ to $T$ defines an isomorphism
$R(G)\otimes\C\to(R(T)\otimes\C)^{W_G}$. Moreover, if $G$ is simply connected,
Theorem 7.6 yields that $R(G)\otimes\C$ is a polynomial algebra in $k=\rank G$
free generators. From representation theory one actually knows
$R(G)\cong\Z[\chi_{\varrho_1},\dots,\chi_{\varrho_k}]$, where
$\varrho_j,j=1,\dots,k$, is the representation with highest weight
$\bar\omega_j$, and that
$\chi_{\varrho_j}\circ\pi_{\vert_{\hat\Sigma}}=S(e^{2\pi i\bar\omega_j})$ +
lower order terms.

\bigskip
The discussion of the last example shows that isoparametric submanifolds are
intimately related to the representation theory of compact Lie groups. This
becomes even more evident by Corollary 7.7. In fact, if $M$ is a compact
isoparametric submanifold with flat sections of a compact simply connected
normal homogeneous space $G/H$, then $M$ is the level set of a function
$\psi=(\psi_1,\dots,\psi_k)$ whose components $\psi_i$ are eigenfunctions of
the Laplacian. The functions $\tilde\psi_i:=\psi_i\circ\pi_G$, where
$\pi_G:G\to G/H$ is the projection, are then eigenfunctions of the
Laplacian of $G$ and hence coefficients of representations of $G$.

\vfill\eject

   \centerline {\bf Appendix A}

\bigskip

   \centerline{\bf A remark on polar actions}

\bigskip
In this appendix we prove a result which has been conjectured by  Palais and Terng in
[PT1], Remark 3.3 and  [PT2], 5.6.8.

\bigskip\noindent
{\bf Theorem A} {\sl Let $X$ be a complete Riemannian manifold  and $G$ a Lie group
which acts properly  by isometries on $X$. If the distribution of normal spaces to the regular
orbits is integrable then there exists a complete totally geodesic immersed section for the
action of $G$ which moreover meets all orbits and always perpendicularly.}

\bigskip

We will  recall some basic facts about
existence and extendibility of totally geodesic submanifolds first and begin with the  following
 result of Cartan.

\bigskip
\noindent{\bf Theorem (E. Cartan, see [Her])}
{\sl Let $M$ be a Riemannian manifold, $p\in M$ and $W$ a linear subspace of
$T_pM$. Then there exists a totally geodesic submanifold $N$ of $M$ with $p \in N$
and $T_pN = W$ if and only if there exists some $\epsilon \in \R_+$
such that for every geodesic $\gamma$ in $M$ with $\gamma(0) = p$
and $\gamma'(0) \in W \cap U_\epsilon(0)$ the Riemannian curvature
tensor of $M$ at $\gamma(1)$ preserves the parallel translate of $W$ along
$\gamma$ from $p$ to $\gamma(1)$.}

\bigskip\noindent
 Let $M$ be a fixed Riemannian manifold. We recall that two isometric immersions
$f_i: N_i \rightarrow M$, $i=1,\, 2$ are said to be equivalent, if there exists a global isometry
$g: N_1\rightarrow N_2$ such that $f_1= f_2 \circ g$.
If $f: N^k \rightarrow M$ is an isometric totally geodesic immersion, then it induces a
 (differentiable) map $\tilde {f}: N \rightarrow G_k(M)$, where
$G_k(M)$ is the Grassmannian of $k-$planes of $TM$.
Namely, $\tilde f (x) = (f(x), f_*(T_xN))$.
We say that the isometric totally geodesic immersion  $(N,f)$ is {\it  compatible}, if $N$ is
 connected and $\tilde {f}$ is injective. Any isometric totally geodesic immersion  from a
 connected Riemannian manifold $N$ can be factorized through a compatible one. In fact, this
 can be done by identifying points in $N$ with the same image in the Grassmannian. In this way the
 quotient space is a differentiable manifold which admits a unique Riemannian structure such
 that the projection from $N$ is a local isometry. This is due to the fact that locally a
 totally geodesic submanifold is completely determined by its tangent space at a point. By
  the same reason one can show the following: a compatible isometric totally geodesic immersion
   is completely determined, up to equivalence, by its image into the Grassmannian.

Let $\T$ denote the collection of all (equivalence classes of) compatible isometric totally
 geodesic immersion  into $M$. $\T$ has a natural partial order $\preceq$. Namely,
 $(N_1, f_1)\preceq
(N_2,f_2)$, if there is a 1-1 local isometry $i: N_1\rightarrow N_2$ such that $f_1 = f_2\circ i$.
 In this case we say that $(N_2,f_2)$
{\it extends} $(N_1,f_1)$. It is standard to show  that $(N_1, f_1)\preceq (N_2,f_2)$, if and only
 if $\tilde f_1(N_1) \subset \tilde f_2(N_2)$. Using Zorn's Lemma, there
 exists for each
 $(N,f)\in \T$, a maximal $(\bar N, \bar f) \in \T$ extending $(N,f)$.

\bigskip\noindent
{\bf Lemma 1} {\sl Let  $(N_1, f_1),(N_2,f_2)\in \T$ and assume that
$\tilde f_1(N_1) \cap \tilde f_2(N_2)\neq \emptyset $. Then there exists $(N,f)\in \T$
which extends both immersions.}

\bigskip\noindent
{\bf Proof}: Define $N$ to be the disjoint union of $N_1$ and $N_2$ with the following relation:
$x_1 \sim \  x_2$ iff $\tilde f _1 (x_1) = \tilde f _2 (x_2)$, for $x_1 \in N_1$ and $x_2 \in N_2$.
 Observe that $N$ is connected. Let $f:N\rightarrow M$ be defined by $f([x_i]) = f_i (x_i)$,
 $x_i\in N_i$, $i=1,2$. The proof follows now in a standard way.\hfill\qed

\bigskip\noindent
{\bf Lemma 2} {\sl Let $(N_1, f_1), (N_2,f_2)\in \T$ with $(N_1,f_1)$ maximal
and assume that there exist sequences
$\{x_n\}$ in $N_1$ and $\{y_n\}$ in $N_2$
such that :

(i) $\tilde f_1 (x_n) = \tilde f_2 (y_n)$.

(ii) $\{ y_n\}$ converges to some $y \in N_2$.

\noindent
Then also $x_n$ converges.}

\bigskip\noindent
{\bf Proof}: Immediate from Lemma 1.\hfill\qed

\bigskip\noindent
{\bf Lemma 3} {\sl Under the assumptions of the above theorem, let $p\in X$ and let $H$ be the
connected component of the isotropy subgroup of $G$ at $p$. Then the isotropy representation of
$H$ on the normal space $\nu_p(G.p)$ (i.e. the so called slice representation)
 is polar. Moreover, if $\Sigma$ is a (linear) section for this representation then
$exp_p (\Sigma \cap B_\varepsilon (0))$ is locally  a (totally geodesic) section for the  action
of $G$, for $\varepsilon$ small.}

\bigskip\noindent
{\bf Proof}: First observe that a totally geodesic connected submanifold
$\tilde \Sigma$ of $X$ meets orbits orthogonally if and only if, for some $q\in \tilde \Sigma$
and any $Z\in Lie (G)$ (identified with a Killing field on $X$), the following two conditions
 hold:

(i) $Z.q \in (T_q\tilde \Sigma)^\perp$.

(ii) $(\nabla Z)_q (T_q\tilde \Sigma) \subset (T_q\tilde \Sigma)^\perp$.

In fact, if $\gamma (t)$ is a geodesic in $\tilde \Sigma$ starting at $q$ then the Jacobi field
$J(t)= Z.\gamma (t)$ has initial conditions $J(0) = Z.q$, $J'(0)= \nabla _{\gamma '(0)} Z$ which
are both in $(T_q\tilde \Sigma)^\perp$.

Choose a sequence $\{p_n\}$, in  the regular points of $M$  for the action of $G$, such that
$p_n \rightarrow p$ and let $\tilde \Sigma _n$ be local sections for the action of $G$ with
$p_n \in \tilde \Sigma$. By choosing, eventually, a subsequence we may assume that
$T_{p_n}\tilde \Sigma = \nu _{p_n}$ converges to some subspace $W$ of $T_pX$, which must be
orthogonal to $T_p(G.p)$. By continuity we obtain, in particular, that
$(\nabla Z)_q (W) \subset (W)^\perp$ for all $Z \in
 Lie (H)$. Since $ {d\over dt}|_0 (exp tZ)_{*p} = (\nabla Z)_p$, we have that the Lie algebra
 $\h $ of the image of the isotropy representation of $H$ at $p$ coincides with
 $\{(\nabla Z)_p: Z \in Lie (H)\}$. So,
$\h (W)\subset W^\perp$, which implies the polarity of the slice representation
(since the codimension of the action of $G$ is the same as the codimension of the
slice
representation). Let now $\Sigma\subset T_pX$ be any section of the isotropy action. We want
to show that the image by  $exp_p$, of some open neighborhood of $0\in \Sigma$, is a totally
geodesic submanifold of $X$. We will need the following auxiliary result which is well known
and standard to prove.

\bigskip\noindent
 {\bf Sublemma} {\sl There exists $\varepsilon > 0$ such that the
codimension of $G. exp _p(v)$ in $X$
 is equal to the codimension of
$H.v$ in $\nu_p(G.P)$, for all
$v\in \Sigma $ with
$||v||< \varepsilon$. In particular, if $v$ is a principal vector for the
slice representation with
$\Vert v\Vert<\epsilon$, then the geodesic $\gamma_v(t)$, $0\neq t\in[-1,1]$,
 consists of principal points for the action of $G$ on $X$.}

\bigskip\noindent
{\bf Proof of Lemma 3 (continued)}:
Let $\varepsilon > 0$ be given by the above sublemma which we may assume smaller than the
injectivity radius at $p$. Let $v\in \Sigma$ be a principal vector for the isotropy
representation with $||v||< \varepsilon$ and consider the geodesic $\gamma _v$. Observe
that $\gamma '(t) \in \nu _{\gamma (t)}$ for all $t\in (0,1]$. In fact, if $Z\in Lie(G)$
then
$\langle Z. \gamma _v(0), \gamma _v' (0)\rangle = 0$ and
$ {d\over dt}\vert_0\langle Z. \gamma _v(t), \gamma _v' (t)\rangle =
\langle \nabla _{\gamma_v'(0)}Z,\gamma_v'(0)\rangle = 0$, by the Killing equation.

We have that $\nu ^t:= \nu _{\gamma _v(t)}$, $0<t\leq 1$, is parallel along $\gamma _v(t)$,
since $\nu$ is a totally geodesic distribution. It is standard to show that the distance
$d( \nu^t, T_{\gamma _v(t)}(exp _p (\Sigma \cap B_\varepsilon (0))))$ tends to zero, if
$t\rightarrow 0$ (regarding, for instance, both subspaces as elements of $\Lambda ^k
(T_{\gamma _v(t)}X)$, where $k = dim (\nu)$). Roughly speaking  $T_{\gamma _v(t)}(exp _p
(\Sigma \cap B_\varepsilon (0)))$ approximates the normal space $\nu ^t$ of the orbit $G.
\gamma _v(t)$, for $t$ small, $t\neq 0$. So, if we set $\nu ^0 = \Sigma$, then $\nu ^t$,
$t\in [0,1]$, is a parallel distribution along $\gamma _v$. Since $\nu$ is totally geodesic,
 we have in addition that $\nu^t$ is invariant under the curvature tensor $R$ of $X$, for
 $t\neq 0$ and by continuity for all $t\in [0,1]$.
Since principal vectors are dense in $\Sigma$ we obtain, by a continuity argument, that for
any $v\in \Sigma$ with $||v||<\varepsilon $ the parallel transport of $\Sigma$ along the
geodesic $\gamma _{v|[0,1]}$ is invariant under the curvature tensor $R$. Then $exp _p(
\Sigma \cap B_\varepsilon (0))$ is a totally geodesic submanifold of $X$ due to the Theorem of
Cartan.
Observe that the parallel transport of $\Sigma$ along the geodesic
$\gamma _{v|[0,1]}$, $v\in \Sigma$ principal, coincides with
$T_{\gamma _v(1)}(exp _p (\Sigma \cap B_\varepsilon (0)))$. By construction, this parallel
transport must also coincide with $\nu _{\gamma _v(1)}$. So,
$exp _p (\Sigma \cap B_\varepsilon (0))$ is a local section.\hfill\qed

\bigskip\noindent
 {\bf Proof of Theorem A}: Let $f: N\rightarrow X$ be a maximal isometric totally geodesic
immersion which extends some local section. Observe that $(N,f)$ meets orbits perpendicularly,
since it does so in an open non empty subset (using the fact that Killing fields are Jacobi along
 geodesics). If $(N,f)$ is not complete there exists a geodesic $\gamma :[0,1)\rightarrow N$
 and
a sequence $\{t_n\}$ in $[0,1)$ which tends to $1$ and such that $\{\gamma (t_n)\}$
  is not convergent in $N$. We have that $f_*(T_{\gamma (t_n)}N)$ converges to $W$, the parallel
  transport along $\beta$ of $f_*(T_{\gamma (0)}N)$, where $\beta :[0,1]\rightarrow X$ is the
   geodesic extending $f\circ \gamma$. By Lemma 3 and its proof we have that  $W$ is a section of
   the isotropy representation at $\beta (1)$. Moreover, $exp _{\beta (1)}(W\cap B_\varepsilon
   (0))$
    is a totally geodesic submanifold of $X$, if $\varepsilon $ is small. If $n$ is large, one has
    that
$(exp _{\beta (1)})_{*v_n}(T_{v_n}W)= f_{*\gamma (t_n)}(T_{\gamma (t_n)}N)$, where $v_n = (1-t_n)
(-\beta ' (1))$ (since both subspaces are parallel along $\beta$). Then one has, for the induced
maps on the Grassmannian, that $\tilde f (\gamma (t_n)) =
\tilde {exp}_{\beta (1)} (v_n)$. Since $v_n$ converges to $0$ and
$\gamma (t_n)$ does not converge, we conclude, by Lemma 2, that $(N, f)$ is not maximal. A
contradiction. Thus $N$ is complete. Observe that $f(N)$ meets every orbit, as $G-$orbits are
 closed (the action is proper) and so, the exponential of the normal space to any orbit meets all
  other orbits.\hfill\qed

\bigskip\noindent
{\bf Remark}: After finishing this paper we learned that K. Grove and W.
Ziller gave an independent proof of Theorem A. This will appear in [GZ].
Moreover G. Thorbergsson draw our attention to a paper by H. Boualem [Bo]
which treats the more general case of singular riemannian foliations and
contains Theorem A as a special case. But we hope that our shorter proof in
the particular situation will still be of some interest.

\vskip 1 cm

\centerline{\bf Appendix B}

\bigskip
{\bf Isoparametric submanifolds in Hilbert spaces are embedded and have globally
flat normal bundle}

\bigskip
Terng has shown that isoparametric submanifolds of euclidean space have
globally flat normal bundle ([T1], Proposition 3.6).
The aim of this appendix is to generalize this result
 to infinite dimensions by a slight modification of her arguments.
Moreover we fill a gap in the proof of Corollary 7.2 of [T2] which says that
isoparametric submanifolds of Hilbert spaces are embedded.
Both problems are related, and will be resolved simultanously. Thus we are
going to prove:

\bigskip\noindent
{\bf Theorem B} {\sl Let $V$ be a separable Hilbert space and $f:M\to V$ a proper
Fredholm immersion of finite codimension which satisfies
\item{(i)} $\nu M$ is flat, and
\item{(ii)} for any parallel normal vector field $\xi (t)$ along any
differentiable curve $c:[0,1]\to M$, the shape operators $A_{\xi(0)}$ and
$A_{\xi(1)}$ are orthogonally equivalent.

\noindent
If $M$ is connected, then $f(M)$ is an embedded submanifold of $V$ with
globally flat normal bundle. In particular $M$ (more precisely $f$) itself has globally flat normal
bundle and $M$  is isoparametric according to the definition of Terng [T2].}

\bigskip\noindent
{\bf Remark}: One has to be a little bit careful with the notion of globally
flat normal bundle of an immersion $f:M\to X$. Even if $f(M)$ is embedded this
depends strongly on $f$. For example, if $\nu M$ has finite holonomy (as in
our case, see the first lines of the proof of Theorem B) one can replace $M$ by a finite
 cover such that the new
immersion has globally flat normal bundle. This shows in particular that the
assumption (i) above may be replaced by global flatness of $\nu M$.

\bigskip
We are going to deduce Theorem B from the following result of Terng which
seems to us the correct statement of what she proved in [T2] about the
embeddedness of isoparametric submanifolds. Namely, to show that an immersed
isoparametric submani\-fold $f:M\to V$ of a Hilbert space is embedded, Terng
proves that each critical point of the euclidean distance function (to a
non-focal point) has a Bott-Samelson cycle and that therefore $f$ ist taut
(Theorem 7.1 of [T2]) and hence an embedding (Corollary 5.9 of [T2]).
However, the construction of the Bott-Samelson cycles only works if $f$ is
injective on all curvature spheres (cf. the proof of Lemma 8.3.2 in [PT2]).
While the injectivity is clear for curvature spheres of dimension at least $2$
(as they are simply connected and the restriction of $f$ is a local isometry
onto a round sphere), it is in general not true for $1$-dimensional curvature
circles.

\bigskip\noindent
{\bf Theorem (Terng [T2])} {\sl Let $V$ be a separable Hilbert space and $f:M\to V$
a proper Fredholm immersion of finite codimension which is isoparametric in
the sense of Terng, i.e. satisfies
\item{(i)} $\nu M$ is globally flat, and
\item{(ii)} for any parallel normal vector field $\xi$ on $M$, the shape
operators $A_{\xi(p)}$ and $A_{\xi(q)}$ are orthogonally equivalent for all $p,q\in M$.

\noindent
If $M$ is connected and $f$ is injective on all curvature circles, then $f$ is an embedding.}

\bigskip\noindent
{\bf Remarks}:

\noindent
(i) We are grateful to G. Thorbergsson who pointed out to us another proof of
this theorem  using only the mountain pass lemma, ([Th]).

\noindent
(ii) If all multiplicities of $M$ (i.e. the dimensions of the curvature
spheres) are at least $2$, then Theorem B is a direct consequence of Terng's result
and the remark following Theorem B.

\noindent
(iii) Under the assumptions of Terng's theorem the parallel immersions
$f_\xi:M\to V,p\mapsto f(p)+\xi(p)$, are also embeddings for all non-focal
$\xi\in\nu (M)$. In fact, this new immersion has also globally flat normal bundle
and is injective on all curvature spheres as well (since along a curvature
sphere $\xi(p)=\lambda(f(p)-m)+\xi_0$ for some $\lambda\in\R$, and $m,\xi_0\in
V$ and hence $f_{\xi}(p)=(1+\lambda)f(p)-\lambda m+\xi_0)$.

\bigskip\noindent
{\bf Proof of Theorem B}: We may assume that $M$ is full and that therefore
$\nu_xM$ is spanned for any $x\in M$ by the
curvature normals $n_i(x),i\in I$. These are a well defined set of normal
vectors at each point, invariant under parallel translation, but not necessarily of
globally defined vector fields. Let $x_0\in M$ and
 $S:=\{n_i(x_0)/\Vert n_i(x_0)\Vert\ge r\}$ for some $r>0$. Since $M$ is full and of finite
codimension we may assume that $r$ is big enough, so that $S$ spans
$\nu_{x_0}(M).\ S$ is a finite set, because the shape operators are compact.
Thus the holonomy group of $\nu M$ at $x_0$, which preserves $S$, is finite.
Moreover it is trivial if the curvature normals have pairwise different
length.

To get theorem B from Terng's result (by going to a finite cover of $M$ we
easily could assume that $\nu M$ is globally flat) we want to factorize $f$ over
an appropriate quotient but run then into the difficulty that the normal
bundle of the quotient might not be globally flat any more.
We therefore assume first that the curvature normals $n_i$ have different
lengths for all $i\in I$. Hence $\nu
(M)$ is globally flat and the $n_i$ are globally defined.
We fix $i\in I$ and let $S_i(q)$ denote the curvature sphere through $q\in M$,
which is by definition the integral manifold of the curvature distribution
$E_i$. Then $f(S_i(q))$ is a round sphere of radius
$1\over\Vert n_i\Vert$ for all $q\in M$.  We  assume $\dim S_i(q)=1$ in
addition and define an equivalence relation $\sim$ on $M$ by $x\sim y$ iff
$f(x)=f(y)$ and $y\in S_i(x)$. Let $k(x)$ denote the number of points
equivalent to $x$. Thus the closed geodesic $S_i(x)$ covers the circle
$f(S_i(x))\quad k(x)$-times under $f$. We fix $x\in M$ and orient $S_i(x)$ by
choosing a tangent vector field $Z$ of constant length $2\pi/\Vert n_i\Vert$.
This vector field can be extended to a vector field $\tilde Z$ on a
neighborhood $U$ of $S_i(x)$ consisting of integral manifolds of $E_i$ with
$\Vert\tilde Z_i\Vert\equiv 2\pi/\Vert n_i\Vert$ and $\tilde Z(y)\in E_i(y)$
for all $y\in U$. Let $\varphi^t:U\to U,t\in \R$, be the associated flow and
$\varphi:=\varphi^1$. Then $f(\varphi^t(y))$ covers the circle $f(S_i(y))$
just once if $t$ runs from $0$ to $1$. In particular $f(\varphi(y))=f(y)$ for
all $y\in U$ and hence $f(\varphi^{k(x)}y)=f(y)$. Since $\varphi^{k(x)}(x)=x$
and $f$ is locally injective there exists a neighborhood $U'$ of $x$ with
$\varphi^{k(x)}(y)=y$ and hence with $k(y)\le k(x)$ for all $y\in U'$. By
shrinking $U'$ further we then find also a neighborhood $U''$ with $k(y)=k(x)$
for all $y\in U''$, as otherwise there would exist a sequence $y_n$
converging to $x$ with $k(y_n)=l<k(x)$ implying $f(\varphi^lx)=f(x)$ which is
a contradiction. Thus $k(x)$ is locally constant and hence constant on $M$ by
connectedness. This implies that $M/\sim$ is locally given as the
quotient of $M$ by the free action of the group $\{\varphi^l\mid
l=1,\dots,k\}$ of diffeomorphism and is hence in a natural way a manifold
 such that $M\to M/\sim$ becomes a $k$-fold
covering (cf. [Bou2], ch. 5.9.5). Furthermore $f$ factorizes through $M/\sim$. Since the fibres of $f$
are finite, there are only finitely many indices $i\in I$, such that $S_i$ is a
curvature circle of $M$ with more than one equivalent point. Thus repeating
the above process finitely many times one gets a quotient manifold $\bar M$ and
an induced isoparametric immersion
 $\bar f:\bar M\to V$ whose normal bundle is still globally
flat (as the curvature normals have different lengths) and which is injective
on all curvature circles. Thus $\bar f$ is an embedding by Terng's theorem and $f(M)=\bar f(\bar
M)$ is an embedded submanifold with globally flat normal bundle.

\bigskip\noindent
In the general situation (where the $\Vert n_i\Vert$ are not necessarily
pairwise different) we choose $q\in M$ and $\xi\in\nu_qM$ such that
$<\xi,n_i(q)>\ne 1$ and
$$
\Vert {n_i(q)\over 1-<\xi,n_i(q)>}\Vert\ne\Vert {n_j(q)\over
1-<\xi,n_j(q)>}\Vert
$$
for all $i, j\in I$ and $j\ne i$. This can be done since $I$ is countable. Let
$M_\xi:=\{(p,\tilde\xi)\in M\times V\mid\tilde\xi\in \nu_pM$ is obtained by
parallel translating $\xi$ in $\nu M$ along any curve in $M$ connecting $q$
with $p\}$ and $f_\xi:M_\xi\to V$ the mapping defined by
$f_\xi(p,\tilde\xi)=f(p)+\tilde\xi$. Then $f_\xi$ is a proper Fredholm
immersion of finite codimension which satisfies (i) and (ii) of the theorem.
Moreover the curvature normals of $M_\xi$ at $(q,\xi)$ are exactly
$$
{n_i(q)\over 1-<\xi,n_i(q)>}\ ,\ i\in I\ ,
$$
and therefore have different lengths. Thus $f_\xi(M_\xi)$ is an embedded
submanifold of $V$ with globally flat normal bundle by the first step. The vector
$-\xi\in\nu_{(q,\xi)}M_\xi$ defines a global parallel normal field along
$M_\xi$ and hence also along $f_\xi(M_\xi)$. The corresponding parallel
manifold is $f(M)$. By the remark (iii) following Terng's theorem, $f(M)$ is
therefore an embedded submanifold with globally flat normal bundle.\hfill\qed

\vskip 1 cm

\centerline{{\bf References}}

\bigskip\noindent
[BC] Bishop, R.L. and Crittenden, R.J., {\sl Geometry of manifolds,} Academic
Press, New York, 1964

\noindent
[BD] Br\\ocker, Th. and tom Dieck, T., {\sl Representations of compact Lie
groups,} Springer 1995

\noindent
[Bo] Boualem, H., {\sl Feuilletages riemanniennes singuliers transversalement
integrables,} Compositio Math. {\bf 95}, 101 - 125

\noindent
[Bou1] Bourbaki, N., {\sl Groupes et alg\`ebres de Lie,} Chapitre 4,5 et 6,
Hermann, 1968

\noindent
[Bou2] Bourbaki, N., {\sl Vari\'et\'es diff\'erentielles et analytiques,}
Hermann 1967

\noindent
[Ca] Cartan, E., {\sl Familles de surfaces isoparam\'etriques dans les espaces
\`a courbure constante,} Annali di Mat. {\bf 17} (1938), 177 - 191

\noindent
[CW] Carter, S. and West, A., {\sl Isoparametric systems and transnormality,}
Proc. London Math. Soc., {\bf 51} (1985), 520 - 542

\noindent
[DS] Dunford, N. and Schwarz, J.T., {\sl Linear Operators,} Part II,
Interscience Publishers Inc. New York, 1958

\noindent
[GW1] Gromoll, D. and Walschap, G., {\sl Metric fibrations in euclidean
space,} Asian J. Math. {\bf 1} (1997), 716 - 728

\noindent
[GW2] Gromoll, D. and Walschap, G., {\sl The metric fibrations in euclidean
space,} preprint

\noindent
[GZ] Grove, K. and Ziller, W., {\sl Characterizations of Polar actions,} in
preparation

\noindent
[Ha] Harle, C.E., {\sl Isoparametric families of submanifolds,} Bol. Soc.
Bras. Math. {\bf 13} (1982), 35 - 48

\noindent
[Hel1] Helgason, S., {\sl Invariants and fundamental functions,} Acta Math.
{\bf 109} (1963), 241 - 258

\noindent
[Hel2] Helgason, S., {\sl Groups and Geometric Analysis}, Academic Press, 1984

\noindent
[Her]   Hermann, R., {\sl Existence in the large of totally geodesic
submanifolds of Riemannian spaces}, Bull. Amer.Math. Soc. {\bf 66}
 (1960), 59-61

\noindent
[HPTT] Heintze, E., Palais, R.S., Terng, C.L. and Thorbergsson, G., {\sl Hyperpolar
Actions on Symmetric Spaces,} Geometry, Topology \& Physics for Raoul Bott, ed.
by S.T. Yau, (1995), 214 - 245

\noindent
[K] Kato, T., {\sl Perturbation theory for linear operators,} Springer 1966

\noindent
[Ko] Kollross, A., {\sl A classification of hyperpolar and cohomogeneity one
actions}, PhD thesis, Augsburg 1998

\noindent
[KT] King, C. and Terng, C.L., {\sl Minimal submanifolds in path space,}
Global Analysis in Modern Mathematics, ed. by K. Uhlenbeck, Publisch or
Perish, 1993, 253 - 282

\noindent
[M] Milnor, J., {\sl Morse Theory,} Princeton University Press, Princeton,
1963

\noindent
[O'N1] O'Neill, B., {\sl The fundamental equations of submersions,} Michigan
Math. J. {\bf 13} (1966), 459 - 469

\noindent
[O'N2] O'Neill, B., {\sl Submersions and geodesics,} Duke Math. J. {\bf 34}
(1967), 363 - 373

\noindent
[PT1] Palais, R.S. and Terng, C.L. {\sl A general theory of canonical forms,}
Trans. Amer. Math. Soc. {\bf 300} (1987), 771 - 789

\noindent
[PT2] Palais, R.  and  Terng, C.-L., {\sl Critical Point Theory and Submanifold geometry,}
Lect. Notes in Math. {\bf 1353},  Berlin Heidelberg New York: Springer 1988

\noindent
[PTh] Podest\`a, F. and Thorbergsson, G., {\sl Polar actions on rank one
symmetric spaces}, preprint

\noindent
[T1] Terng, C.L., {\sl Isoparametric Submanifolds and their Coxeter Groups,}
J. Differential Geometry {\bf 21} (1985), 79 - 107

\noindent
[T2] Terng, C.L., {\sl Proper Fredholm Submanifolds of Hilbert Space,} J.
Differential Geo\-metry {\bf 29} (1989), 9 - 47

\noindent
[T3] Terng, C.L., {\sl Polar actions on Hilbert space,} J. Geom. Anal. {\bf
5} (1995), 129 - 150

\noindent
[Th] Thorbergsson, G., {\sl private communication}, May 1999

\noindent
[TT] Terng, C.L. and Thorbergsson, G., {\sl Submanifold geometry in symmetric
spaces,} J. Differential Geometry {\bf 42} (1995), 665 - 718

\noindent
[W] Wang, Q., {\sl Isoparametric hypersurfaces in complex projective spaces},
1980 Beijing Symposium on Diff. Geom. and Diff. Eq., Vol. {\bf 3}, 1509 -
1524, Science Press, Beijing 1982

\vskip 2 cm

{\settabs 2 \columns
\+Ernst Heintze                      &\qquad Xiaobo Liu\cr
\+Institut f\\ur Mathematik          &\qquad Department of Mathematics\cr
\+Universit\\at Augsburg             &\qquad University of Notre Dame\cr
\+Universit\\atsstrasse 14           &\qquad Notre Dame, IN 46556-0398\cr
\+D - 86159 Augsburg, Germany        &\qquad USA\cr
\+{\sl heintze@math.uni-augsburg.de} &\qquad{\sl
xliu3@hilbert.helios.nd.edu}\cr}

\vskip 1 cm

{\obeylines\parindent=0pt
\hskip 4 cm Carlos Olmos
\hskip 4 cm Facultad de Matem\'atica, Astronom\'\i a y F\'\i sica
\hskip 4 cm Universidad Nacional C\'ordoba
\hskip 4 cm Medina Allende y Haya de la Torre
\hskip 4 cm Ciudad universitaria
\hskip 4 cm 5000. C\'ordoba, Argentina
\hskip 4 cm {\sl olmos@mate.uncor.edu}
}

\vfill\bye